\documentclass{gtart}
\usepackage{amssymb}

\usepackage{graphicx}
\usepackage{amsmath}
\usepackage{amsthm}
\usepackage{color}
\usepackage{pictex}

\theoremstyle{plain} 

\newtheorem{theorem}{Theorem}[section]
\newtheorem{thm}[theorem]{Theorem}
\newtheorem{lemma}[theorem]{Lemma}

\newtheorem{prop}[theorem]{Proposition}

\theoremstyle{definition}
\newtheorem{defn}[theorem]{Definition}

\newtheorem{example}[theorem]{Example}
\newtheorem{question}[theorem]{Question}

\theoremstyle{remark}
\newtheorem{rmk}[theorem]{Remark}

\newtheorem{obs}[theorem]{Observation}

\newcommand{\be}{\begin{equation}}
\newcommand{\ee}{\end{equation}}
\newcommand{\ba}{\begin{align}}
\newcommand{\ea}{\end{align}}
\newcommand{\ben}{\begin{enumerate}}
\newcommand{\een}{\end{enumerate}}
\newcommand{\bi}{\begin{itemize}}
\newcommand{\ei}{\end{itemize}}
\newcommand{\ra}{{\rightarrow\,}}

\providecommand{\abs}[1]{\lvert#1\rvert}

\providecommand{\gint}[1]{\lfloor#1\rfloor}
\providecommand{\lint}[1]{\lceil#1\rceil}

\newcommand{\bC}{{\mathbb C}}

\newcommand{\bN}{{\mathbb N}}
\newcommand{\bR}{{\mathbb R}}
\newcommand{\bZ}{{\mathbb Z}}
\newcommand{\bQ}{{\mathbb Q}}
\newcommand{\cA}{{\cal A}}
\newcommand{\cC}{{\cal C}}
\newcommand{\cG}{{\cal G}}
\newcommand{\cL}{{\cal L}}

\newcommand{\cS}{{\cal S}}
\newcommand{\cT}{{\cal T}}
\newcommand{\bdy}{\partial}
\newcommand{\del}{\partial}
\newcommand{\dponedX}{{\frac{\del p_1}{\del X}}}
\newcommand{\dponedY}{{\frac{\del p_1}{\del Y}}}
\newcommand{\dptwodX}{{\frac{\del p_2}{\del X}}}
\newcommand{\dptwodY}{{\frac{\del p_2}{\del Y}}}

\begin{document}

\title{Toric structures on near-sym\-plec\-tic 4-manifolds}
\shorttitle{}

\author{ David T. Gay
\footnote{Supported in part by NSF/DMS-0244558 and CRM/ISM}
\\
    Margaret Symington
\footnote{Supported in part by NSF/DMS-0204368, AIM, and NSF FRG} 
}
\address{Department of Mathematics and Applied Mathematics,
University of Cape Town,
Private Bag X3,
Rondebosch 7701,
South Africa}
\secondaddress{Department of Mathematics,
Mercer University,
1400 Coleman Ave.,
Macon, GA 31207, USA}
\email{dgay@maths.uct.ac.za}
\secondemail{symington\_mf@mercer.edu}

\begin{abstract}
A near-symplectic structure on a 4-manifold is a closed 2-form that is
  symplectic away from the 1-dimensional submanifold along which it
  vanishes and that satisfies a certain transversality condition along
  this vanishing locus.  We investigate near-symplectic 4-manifolds
  equipped with singular Lagrangian torus fibrations which are locally
  induced by effective Hamiltonian torus actions. We show how such a
  structure is completely characterized by a singular integral affine
  structure on the base of the fibration whenever the vanishing locus
  is nonempty.  The base equipped with this geometric structure
  generalizes the moment map image of a toric 4-manifold in the spirit
  of earlier work by the second author on almost toric symplectic
  4-manifolds. We use the geometric structure on the base to
  investigate the problem of making given smooth torus actions on
  4-manifolds symplectic or Hamiltonian with respect to
  near-symplectic structures and to give interesting constructions of
  structures which are locally given by torus actions but have
  nontrivial global monodromy.
\end{abstract}



\sloppy

\maketitlepage

\section{Introduction}

Advances in symplectic topology in the last decade have shown that
symplectic four-manifolds populate a vast portion of the world of
smooth four-manifolds, extending far beyond the class of K\"ahler
manifolds (see~\cite{fintstern, gompf} among others).
Meanwhile, some of the most powerful techniques for studying them are
motivated by complex algebraic geometry.  For a two-form $\omega$ on a
$2n$-dimensional manifold to be symplectic, it must be non-degenerate
(i.e. $\omega^n$ must be non-vanishing) and it must be closed (i.e.
$d\omega=0$).  The first condition guarantees that the manifold admits
an almost complex structure, while the additional condition of
closedness allows one to get control of solutions to differential
equations involving an \lq\lq approximate Cauchy-Riemann\rq\rq\ 
operator and obtain compact moduli spaces of solutions.

Currently it is becoming apparent that one can study a more general
class of oriented smooth 4-manifolds using the techniques of
pseudo-holomorphic curves~\cite{Taubes} and Lefschetz
fibrations~\cite{ADK}.  Indeed, it suffices that there be a cohomology
class whose square induces the given orientation of the 4-manifold
$X$, i.e $b^+_2(X)>0$.  Honda~\cite{Honda.generic} showed that on such
a manifold there always exists what is now known as a near-symplectic
form:

\begin{defn} \label{D:near-sympl}
On a smooth, oriented $4$-manifold $X$, consider a closed $2$-form $\omega$
such that $\omega^2\ge 0$ and let $Z_\omega$ denote the {\em vanishing
locus}, the set of points where $\omega=0$.
The form $\omega$ is {\em near-symplectic} if
\begin{enumerate}
\item $\omega^2>0$ on the complement of $Z_\omega$ and
\item at each point $x\in Z_\omega$, if we use local coordinates on a
  neighborhood $U$ of $x$ to identify the map $\omega: U \rightarrow
  \Lambda^2 (T^* U)$ as a smooth map $\omega: \bR^4 \rightarrow
  \bR^6$, then its linearization at $x$, $D \omega_x : \bR^4 \ra
  \bR^6$, has rank $3$.
\end{enumerate}
We call the set $X \setminus Z_\omega$ the {\em symplectic locus}.  A
{\em near-symplectic manifold} is an oriented smooth 4-manifold
equipped with a near-symplectic form.
\end{defn}
Definition~\ref{D:near-sympl} is a rephrasing, in local coordinates,
of the definition of a near-symplectic form given in~\cite{ADK}.
Indeed, if $\phi:\bR^4\ra\bR^4$ is a change of coordinates on $U$ and
$\Phi:\bR^6\ra\bR^6$ is the corresponding change of coordinates on
$\Lambda^2(T^*U)$ then, restricted to $Z_\omega$,
$D_x\omega\circ\phi=\Phi\circ D_x\omega$.  Consequently, on $Z_\omega$
we have that $D_x\omega$ represents an intrinsically defined
derivative, the derivative denoted in~\cite{ADK} by $\nabla
\omega_x:T_xX\ra\Lambda^2(T_x^*X)$.  The same paper explains why $3$
is the maximum possible rank for $D \omega_x$ and that this definition
is equivalent to the original definition (see~\cite{Honda.generic}) in
terms of metric properties.  Specifically, a closed $2$-form $\omega$
on a $4$-manifold $X$ is near-symplectic if it is self-dual with
respect to some metric $g$ and, viewed as a section of the bundle
$\Lambda_2^+$ of $g$-self-dual $2$-forms, is transverse to the zero
section.

Several threads of emerging research indicate that one should be able
to understand the moduli spaces of pseudo-holomorphic curves in
near-symplectic manifolds and that the extra structure of a fibration
induced by a Hamiltonian torus action should aid in this endeavor.
Specifically, Taubes has made initial steps in his
program to develop Gromov-Witten invariants for near-sym\-plec\-tic
manifolds that should be invariants of the underlying smooth
structure~\cite{Taubes, Taubes.beasts}; Mikhalkin has calculated, via tropical
algebraic geometry, the Gromov-Witten invariants of toric surfaces in
terms of $1$-complexes in their moment map
images~\cite{Mikhalkin.tropicalenumeration}; and Parker has used
symplectic field theory to gain an understanding of moduli spaces of
pseudo-holomorphic curves in $T^*T^2$ in terms of $1$-complexes in
$\bR^2$~\cite{Parker.1complexes}.  Presuming success on these fronts,
the results in this paper lead one to expect to be able to calculate Taubes'
invariants for locally toric near-sym\-plec\-tic manifolds
by counting $1$-complexes suitably immersed in the bases
of the induced fibrations.

This paper can usefully be read in parallel with the work of
Kaufman~\cite{Kaufman}, which also develops a theory of toric
structures on near-symplectic manifolds.

We begin by giving a characterization of Hamiltonian torus actions on
symplectic 4-manifolds that is
convenient for generalizing to the near-symplectic setting and that
emphasizes the induced fibration.

The following proposition
can be read as a definition by those unfamiliar with Hamiltonian group
actions, but is really a statement that a certain characterization of
\lq\lq Hamiltonian\rq\rq\ is equivalent to the standard definition
(see Definition~\ref{D:Hamiltonian}).

\begin{prop} \label{P:Hamiltonian}
A smooth torus action $\sigma: T^2 \times X \ra X$ on a symplectic
$4$-manifold $(X,\omega)$ is {\em Hamiltonian} if and only if
there exists a smooth map $\mu : X \ra \bR^2$, called the {\em moment map},
such that, for any $\xi \in \bR^2$, the vector field $V_\xi$ whose flow is
$x \mapsto \sigma(t\xi,x)$ (for $t \in \bR$) is defined by
the equation $\omega(V_\xi, W) = -\xi \cdot
D\mu(W)$ for all $W \in TX$.
\end{prop}
A proof of this proposition is provided for convenience in the
appendix.  Note that this proposition is true only for torus actions.
Indeed, the proof relies on the fact that the group is both abelian
and compact.

\begin{rmk}
As part of the proof of Proposition~\ref{P:Hamiltonian}, we establish that
$\omega(V_\xi,V_\eta)=0$ for all $\xi,\eta\in\bR^2$.  This implies that
preimages of the moment map are isotropic, and in particular the top 
dimensional preimages are Lagrangian.
\end{rmk}

Generalizing to the near-symplectic setting we have:
\begin{defn} \label{D:NearSymplecticHamiltonian}
  A smooth torus action $\sigma: T^2 \times X \ra X$ on a
  near-symplectic $4$-manifold $(X,\omega)$ which preserves $\omega$
  is {\em Hamiltonian} if there exists a smooth map $\mu: X \ra \bR^2$
  such that $\mu|_{X \setminus Z_\omega}$ is a moment map for
  $\sigma|_{X \setminus Z_\omega}$.  In particular, $\sigma|_{X
  \setminus Z_\omega}$ is Hamiltonian in the usual sense.  We call
  $\mu: X \ra \bR^2$ the {\em moment map} for $\sigma$.
  An action is {\em locally Hamiltonian} if every orbit has an open
  neighborhood in which the action is Hamiltonian.
\end{defn}

The assumption of the smoothness of $\mu$ is not constraining.
Indeed, given any moment map $\mu: X \setminus Z_\omega \rightarrow \bR^2$
for the restriction of a smooth torus action to the symplectic locus of a
near-symplectic manifold $(X,\omega)$, it can be shown that $\mu$
extends smoothly across $Z_\omega$.

Two questions we address in this paper are:

\begin{question} \label{Q:Hamiltonian}
  Which closed {\em $T^2$-manifolds} ($4$-manifolds equipped with smooth 
effective torus actions) admit near-symplectic forms with respect to
  which the actions are Hamiltonian?
\end{question}

\begin{question} \label{Q:LocallyHamiltonian}
Which closed $T^2$-manifolds admit near-symplectic forms with respect
to which the actions are locally Hamiltonian?
\end{question}

Recall that a {\em toric manifold} is a symplectic manifold equipped
with an effective Hamiltonian torus action (i.e. only the identity
acts trivially).  This definition generalizes immediately to
near-symplectic manifolds.  

We are particularly interested in the
fibrations induced by Hamiltonian and locally Hamiltonian torus actions.
These structures, rather than the actions themselves, appear to be useful for
the study of pseudo-holomorphic curves~\cite{Parker.1complexes}.
Furthermore, shifting the focus to fibrations, there is no need for the
manifold to admit a global effective torus action, thereby extending the
scope of our results.  Therefore we make the following definitions.

\begin{defn}
\label{D:corners}
The boundary of a smooth surface has {\em corners} if its boundary is
piecewise smooth and each nonsmooth point of the boundary
has a neighborhood that smoothly surjects onto a neighborhood of the vertex
of a sector in $\bR^2$.
\end{defn}

\begin{defn} \label{D:toric} Let $(X,\omega)$ be a near-symplectic
  manifold. A {\em toric fibration} of a near-symplectic manifold
  is a smooth surjective map $\pi:X\ra B$ to a smooth surface
  with boundary and corners such that the top-dimensional fibers of
  $\pi$ are Lagrangian tori and, ignoring the smooth structure on $B$, the
  map $\pi$ is the orbit space projection for an effective
  Hamiltonian torus action on $(X,\omega)$.
  A {\em locally toric 
  fibration} of $(X,\omega)$ is a smooth surjective map $\pi:X\ra B$
  to a smooth surface with boundary and corners such that each fiber has an
  open neighborhood in which the fibration is toric for some
  Hamiltonian torus action on the neighborhood.  A {\em (locally) toric
    near-symplectic manifold} is a near-symplectic manifold equipped
  with a (locally) toric fibration.
\end{defn}

These generalizations lead us to the following question:
\begin{question} \label{Q:locallytoric}
  Which smooth 4-manifolds with $b_2^+>0$ can be equipped with
  near-symplectic forms so as to admit locally toric fibrations?
\end{question}

Our central results (Theorems~\ref{T:ToricBases} and~\ref{T:ToricUniqueness}) 
assert that toric near-symplectic structures are in one-to-one
correspondence with surfaces equipped with 
certain singular integral affine structures.  
(See Definition~\ref{D:IASEdgeFolds}).  
This result generalizes naturally to locally toric near-symplectic manifolds
(Theorems~\ref{T:LocallyToricBases} and~\ref{T:LocallyToricUniqueness}).
These theorems allow us to answer
questions about realizing smooth torus actions in the near-symplectic
world and give a number of interesting constructions of
near-symplectic 4-manifolds.

For simply connected manifolds we can answer Questions~\ref{Q:Hamiltonian},
\ref{Q:LocallyHamiltonian} and \ref{Q:locallytoric} concisely:

\begin{thm} \label{T:SimplestThm}
  Every locally toric fibration of a simply connected near-symplectic
  manifold is toric.

  Every smooth effective torus action on a simply connected
  $4$-manifold $X$ with $b_2^+(X)>0$ is Hamiltonian with respect to
  some near-symplectic structure.  Furthermore, the vanishing locus
  for any such near-symplectic structure must have exactly
  $b_2^+(X)-1$ components.
\end{thm}
In~\cite{OrlikRaymond}, which we use extensively here, it is shown
that the simply connected $T^2$-manifolds are precisely $S^4$, 
$S^2\times S^2$, and all connected sums of $\bC P^2$ and
$\overline{\bC P^2}$.   Of course, to support a
near-symplectic structure the manifold must either be $S^2\times S^2$ or have 
at least one $\bC P^2$ summand to make $b_2^+ > 0$.

To state other results requires the notion of positive turning
along boundary components of the orbit space of a torus action
(Definition~\ref{D:PositiveTurning}). Briefly, each boundary component
of the orbit space is naturally decomposed as a union of edges, to each of
which is associated a \lq\lq slope\rq\rq\ in $\bQ \cup \{\infty\}$,
corresponding to the stabilizer subgroups of orbits above the segments
(as shown in~\cite{OrlikRaymond}).  The orbit space, together with this
decomposition, is the {\em weighted orbit space}.  As one traverses a boundary
component, the {\em positive turning} is the total angular turning of these
slopes, where from one edge to the next one turns always
counterclockwise.  Because each boundary component is a closed curve, 
this number is always a non-negative integral multiple of $\pi$.

The following proposition shows how the notion of positive turning
can be used to recognize the topological type of a $T^2$-manifold from its
orbit space.  Accordingly, this result is complementary to those 
in~\cite{OrlikRaymond}.  

\begin{prop}
\label{P:SimplyConnectedT2} 
Consider a simply connected $T^2$-manifold $(X,\sigma)$ whose
weighted orbit space has more than four edges.  Let $T$ be the positive turning
of the weighted 
orbit space and let $V$ be the number of its vertices.  Then $X$ is 
diffeomorphic
to $m\bC P^2\# n\overline{\bC P}^2$ where $m=\frac{1}{\pi}T-1$ and 
$n=V-m-2$.
\end{prop}
Note that when there are four or fewer edges, the topology can be
determined merely by comparing the weighted orbit space with a short
set of examples.

\begin{thm} \label{T:HamiltActions} Consider a closed $T^2$-manifold 
  $(X,\sigma)$ such that   $\sigma$ has no nontrivial finite
  stabilizers.    Let $T_0$ be the largest positive turning along
  any of the boundary components of the orbit space $B$.
  Then there exist near-symplectic forms with respect to
  which $\sigma$ is Hamiltonian only if $T_0\ge2\pi$.  If $g=0$ then
  having $T_0 \geq 2\pi$ is in fact sufficient to give the existence
  of such a near-symplectic form.

  Meanwhile, given any compact surface with nonempty boundary, 
  there is a toric 
  near-symplectic manifold with orbit space $B$.
\end{thm}

\begin{question}
In the preceding theorem, when $g > 0$, is there a clean statement of
a necessary and sufficient condition for making the action Ham\-il\-to\-nian?
\end{question}
The second author is investigating this question.  In the meantime, we
can provide necessary and sufficient conditions for making the action
locally Hamiltonian:

\begin{thm}  \label{T:locallyHamilt}
  Consider a closed $T^2$-manifold $(X,\sigma)$ whose action $\sigma$ has
  no nontrivial finite stabilizers.
  Let $g$ and $k$ be the genus and Euler characteristic of the orbit space
  $B$.
  Let $T$ be the sum of the
  positive turnings along all of the boundary components of the weighted
  orbit space.  Then there
  exist near-symplectic forms on $X$ with respect to which $\sigma$ 
  is locally Hamiltonian if and only if 
\begin{enumerate}
\item $B$ is a torus or
\item $B$ has nonempty boundary and
\ben
\item $g\ge 1$, 
\item $g=0$ and  $T \ge (3-k)\pi$, or 
\item $B$ is an annulus, $T=0$, and the two boundary components each have
one edge of the same slope.
\een
\end{enumerate}
The number of components of the vanishing locus must be 
$\displaystyle\frac{1}{\pi}(T-2\pi\chi)$ where $\chi=2-2g-k$ is the
Euler characteristic of the surface.
\end{thm}

Theorem~\ref{T:locallyHamilt} implies that, up to homeomorphism,
any surface with non-empty boundary is the orbit space of a
near-symplectic manifold with a torus action that is locally Hamiltonian.
Furthermore, there are no constraints on the one-dimensional 
stabilizer subgroups unless $B$ is a disk or an annulus.

A locally toric fibration can arise from a global torus action only if
the integral affine monodromy of the base
(Definition~\ref{D:Monodromy}) is trivial.  With regard to locally
toric fibrations not induced by global torus actions we show:
\begin{prop} \label{P:AnyMonodromy}
  Given any compact surface with boundary $B$
  and any homomorphism $h: \pi_1(B) \ra GL(2,\bZ)$, there is a 
  closed near-symplectic manifold $(X,\omega)$ that admits a 
  locally toric fibration with base $B$ and monodromy
  $h$.  (Here, as part of the construction, $B$ inherits a smooth structure.)
\end{prop}

\begin{prop} \label{P:InfiniteFamily}
  There is an infinite family of mutually non-diffeomorphic closed
  near-symplectic manifolds that support locally toric structures
   with nontrivial monodromy, none
  of which are diffeomorphic to a locally toric near-symplectic manifold 
  with trivial monodromy.
\end{prop}

The general philosophy in this paper, thanks to
Theorems~\ref{T:ToricBases}, \ref{T:ToricUniqueness}, 
\ref{T:LocallyToricBases} and \ref{T:LocallyToricUniqueness}, 
is to study toric and locally toric near-symplectic manifolds 
in terms of immersed polygons in the plane that represent the base, or
a fundamental domain of the base, of the induced fibration.

The bases of toric and locally toric fibrations carry a naturally defined
geometry whose local isometries belong to
${\rm Aff}(2,\bZ) := \{p \mapsto Ap +
b \mid A \in GL(2,\bZ), b \in \bR^2\}$, namely an {\em integral affine
structure}.  Isometric immersions of such surfaces into the plane, 
equipped with
the standard integral affine structure (Definition~\ref{D:IAS}), 
are such that almost every
point has a neighborhood whose image is equivalent, up to the action
of an element of ${\rm Aff}(2,\bZ)$, to a domain in the first quadrant $Q$
of $\bR^2$.  At the remaining points, all of which are on the boundary, 
singularities are allowed
where the boundary can double back on itself as in
Figure~\ref{F:FoldExamples}.  

\begin{example} \label{E:Assorted}
  To conclude this introduction, Figure~\ref{F:FoldExamples} gives an
  indication of the variety of integral affine surfaces with edge
  folds, and hence of locally toric near-symplectic manifolds.  Each
  open circle represents a component of the vanishing locus.  The
  parts of an edge that limit onto an open circle are drawn slightly
  displaced as a visual aid, although they should be understood to
  coincide.  By Theorem~\ref{T:LocallyToricUniqueness} each of these surfaces
  determines a unique locally toric near-symplectic manifold.

\begin{figure}
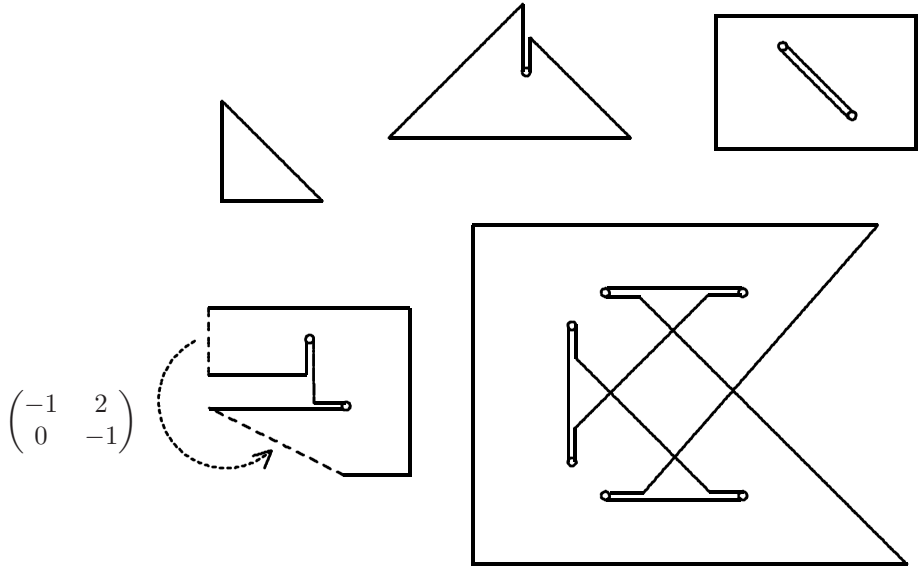

\begin{center}
\include{FoldExamples}
\caption{Various examples of integral affine surfaces, with and
  without edge folds.}
\label{F:FoldExamples}
\end{center}
\end{figure}

  All but one of the figures is the image in $(\bR^2,\cA_0)$ of the base of a
  locally toric fibration via an immersion that preserves integral affine
  structures.  The figure with the dotted lines is the image under such
  an immersion of the complement of a line segment, across which there
  is nontrivial monodromy.  To reconstruct
  the surface, identify the dotted edges by the element of ${\rm Aff}(2,\bZ)$
  comprised of the linear map indicated and an appropriate translation.

  From left to right in the top row, the first figure is the usual
  moment map image of $\bC P^2$ while the second is an integral
  affine base of $\bC P^2 \# \bC P^2$.  (Seeing this would be a good
  exercise to test one's understanding of
  Section~\ref{S:ToricAndIAS}.)  The bottom right figure is the
  immersed image of a surface of genus one with one boundary
  component.  It can be modified easily to give  examples
  with higher genus or more boundary components.
\end{example}

\section{Local models} \label{S:LocalModels}

We devote this section
to models for the neighborhoods of orbits in a near-symplectic manifold
equipped with an effective Hamiltonian torus action. 
There are four types of
orbits that can appear in this setting, characterized by the dimension
of the stabilizer subgroup and whether or not the orbit belongs to the
vanishing locus.

In general, for a smooth torus action on a four-manifold, the possible
stabilizer subgroups are the identity, a circle subgroup, the whole
torus, or a nontrivial finite subgroup~\cite{OrlikRaymond}.
Following~\cite{OrlikRaymond} we denote circle subgroups (stabilizers of
circle orbits) by $G_{(a,b)}:=\{(t_1,t_2) | (t_1,t_2) \cdot (a,b) =
0\}$.  Here and throughout this paper, $(t_1,t_2)$ are 
$\bR/2\pi\bZ$-valued coordinates on the $2$-dimensional Lie group 
$T^2 = S^1 \times S^1$.

An important feature of a Hamiltonian torus action is that the orbits
are {\em isotropic} (i.e., the symplectic form evaluates trivially on
pairs of vectors tangent to an orbit).  Hence, for an effective
Hamiltonian torus action on a four-manifold, the generic orbits are
{\em Lagrangian tori} (isotropic and half the dimension of the ambient
manifold).  Arnold's Theorem~\cite{Arnold} asserts that Lagrangian
torus orbits have a standard product neighborhood, and thereby
prevents the presence of orbits with nontrivial finite stabilizer.
Therefore, in the symplectic locus there are isolated point fibers,
circle fibers that come in one-dimensional families, and the generic
torus orbits.  Meanwhile, because a Hamiltonian action on a
near-symplectic manifold preserves the near-symplectic form, the
vanishing locus must be a union of orbits.  The fact that point orbits
are isolated then implies that each component of the vanishing locus
is one circle orbit.

We will see that the following two examples provide all
the local information about Hamiltonian torus actions on
near-symplectic manifolds.

\begin{example}[Standard torus action on $\bR^4$] \label{E:stdR4}
  Convenient coordinates on $\bR^4 = \bR^2 \times \bR^2$ are the
  square polar coordinates $(p,q):=(p_1,q_1,p_2,q_2)$ which, with
  respect to polar coordinates $(r,\theta)$, are given by
  $p=\frac{1}{2}r^2$, $q=\theta$.  Then the {\em standard symplectic
    structure} on $\bR^4$, which with respect to Euclidean coordinates
  $(x,y)$ is $ \omega_0 = dx\wedge dy := dx_1 \wedge dy_1 + dx_2
  \wedge dy_2$, takes the form $\omega_0=dp\wedge dq := dp_1 \wedge
  dq_1 + dp_2 \wedge dq_2$.  On $(\bR^4,\omega_0)$ we have the {\em
    standard torus action} given by \[t\cdot(p,q)=(p,q+t)\] and the
  corresponding {\em standard moment map} $\mu_0 : \bR^4 \rightarrow Q
  = \{ (x,y) \mid x \geq 0, y \geq 0 \}$ given by $\mu_0(p,q)=p$, so
  $x=p_1$ and $y=p_2$.
\end{example}

The closed first quadrant $Q$ is the orbit space of this action.
For the reader unfamiliar with toric manifolds, we point out that:
\begin{enumerate}
\item The vertex of $Q$ is the image of a point orbit, with torus stabilizer.
\item A non-vertex point on the boundary of $Q$ is the image of a circle
  orbit with stabilizer $G_{(1,0)}$ or $G_{(0,1)}$ depending on
  whether the point is in the positive $x$-axis or positive $y$-axis,
  respectively.
\item A point in the interior of $Q$ is the image of a torus orbit,
  with trivial stabilizer.
\end{enumerate}

\begin{example}[Standard toric action near the vanishing locus]
\label{E:stdvanishinglocus}
Following~\cite{GayKirby.construction} we construct a model
neighborhood of a component of the vanishing locus in a
near-symplectic manifold as follows: Let $\alpha \in S^1$ be the
$2\pi$-periodic coordinate on $S^1$ and let $(x,y,z)$ be Euclidean
coordinates on $\bR^3$.  Then
\[
\omega_1 = 2z (dz\wedge d\alpha + dx\wedge dy) + x (dz\wedge dy -
dx\wedge d\alpha) - y (dz\wedge dx + dy\wedge d\alpha)
\]
is a near-symplectic form on $S^1\times \bR^3$ with vanishing locus
$Z_{\omega_1} = S^1 \times\{(0,0,0)\}$.  Indeed, $\omega_1$ is
self-dual with respect to $g = dz^2 + d\alpha^2 + dx^2 + dy^2$ and
transverse to the zero section of $\Lambda_2^+$.  (Note that the
symplectic orientation is opposite to the standard orientation on $S^1
\times \bR^3$.)

Letting $(r,\theta)$ be polar coordinates in the $(x,y)$-plane in
$\bR^3$, the form $\omega_1$ is invariant under the torus action
$(t_1,t_2) \cdot (\alpha,r,\theta,z) = (\alpha+t_1,r,\theta+t_1,z)$.
The orbit space can be identified with the closed upper half plane $H
= \{ (X,Y) \mid Y \geq 0\}$, so that we have a singular fibration
$\pi: S^1 \times \bR^3 \rightarrow H$ whose fibers are orbits, which
we may take to be given by the equations $X = z$ and $Y = r^2/2)$. This
particular choice of parameterization facilitates certain
calculations in the next section.

Note that if $(X,Y)$ is a point in the interior of $H$ then $\pi^{-1}(X,Y)$ is
a torus orbit with trivial stabilizer, while if $(X,Y)$ lies on the
$X$-axis then $\pi^{-1}(X,Y)$ is a circle orbit with stabilizer
$G_{(1,0)}$.

While the half plane $H$ is the orbit space, it is not the image of
the moment map.  To find the moment map, we find action-angle coordinates --
coordinates on the union of principal orbits such 
that the moment map is projection to the linear coordinates.
(Compare with Example~\ref{E:stdR4}.) 
Accordingly,  define
$f: H \rightarrow \bR^2$ by $p_1=X^2 - Y$, $p_2=2XY$ and let
$\mu_1 = f \circ \pi : S^1 \times \bR^3 \rightarrow \bR^2$.  
Because $\mu_1$ is a submersion onto the complement of the positive $p_1$-axis
and is torus invariant, we can pull $p_1,p_2$ back via $\mu_1$ and get 
coordinates
\[
p_1 = z^2 - \frac{r^2}{2},  \quad q_1 = \alpha, \quad
p_2 = zr^2, \quad q_2 = \theta
\]
on the complement of $S^1\times\bR\times\{0,0\}$, with respect to which
$\omega_1 = dp \wedge dq = dp_1 \wedge dq_1 + dp_2 \wedge dq_2$.
The torus action then becomes $t\cdot(p,q)=(p,q+t)$ and
the image of the union of circle orbits with stabilizer $G_{(1,0)}$ is 
the non-negative $p_1$ axis. 
Consequently, as the reader can verify explicitly, 
$\mu_1$ is the moment map for the action on all of
$S^1 \times (\bR^3 \setminus \{(0,0,0)\})$, and hence a
near-symplectic moment map on $S^1 \times \bR^3$.
\end{example}

\begin{figure}
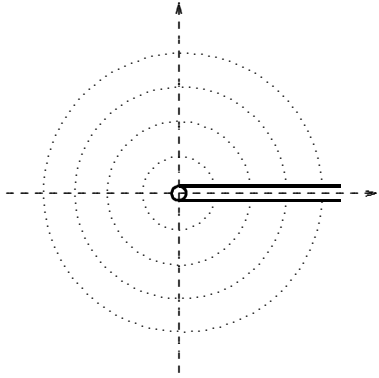

\begin{center}
\include{onefold}
\caption{The moment map image of the standard toric action near a
  component of the vanishing locus in a toric near-symplectic
  manifold.}
\label{F:onefold}
\end{center}
\end{figure}

The preimage under $\mu_1$ of a point on the positive
$p_1$-axis is a disjoint union of two circle orbits and the preimage
of the origin is one circle orbit, the one belonging to
$Z_{\omega_1}$.  To emphasize these features, we draw the moment map
image as the $(p_1,p_2)$-plane with a double line along the positive
$p_1$-axis and a hole at the origin; see Figure~\ref{F:onefold}.

We now show that, up to an automorphism of the torus, 
these examples provide a complete set of examples.

First, it is important to understand the effect of an automorphism
of the torus on the moment map.  
The following lemma is standard for symplectic manifolds (and is
easily verified), and extends to near-symplectic manifolds by continuity.
\begin{lemma} \label{L:AffineMoment} The moment map $\mu$ for a toric
  near-symplectic manifold $(X,\omega,\sigma)$ (where $\sigma$ is the
  torus action) is unique up to addition of a constant $b\in\bR^2$.
  Furthermore, if $t \mapsto A t$, $A\in GL(2,\bZ)$ is any
  automorphism of the torus, then the toric manifold
  $(X,\omega,\sigma')$, with action $\sigma'=\sigma\circ(A\times Id)$,
  has moment map $A^{T}\circ\mu$ where $A^{T}$ is the transpose of $A$.
\end{lemma}

Note that the set of orbit preserving symplectomorphisms is precisely
the set of equivalence classes of equivariant symplectomorphisms in which
two equivariant symplectomorphisms are deemed equivalent if one can
be obtained from the other by precomposing with an automorphism of the
torus.

\begin{lemma} \label{L:stdtoricmodeluniqueness}
  Each orbit in a toric (symplectic) $4$-manifold has a
  torus-invariant neighborhood that symplectically embeds, in an
  orbit-preserving fashion, into $(\bR^4,\omega_0)$ equipped with the
  standard Hamiltonian torus action described in
  Example~\ref{E:stdR4}.  Such an embedding can be chosen to be
  equivariant if and only if the stabilizer subgroup is the identity,
  $G_{(1,0)}$, $G_{(0,1)}$, or the whole torus.
\end{lemma}

\begin{proof}
The equivariant tubular neighborhood theorem (or slice theorem) states
that any orbit has a neighborhood that is equivariantly
symplectomorphic to a neighborhood of the zero section of its normal
bundle equipped with a linear Hamiltonian torus action
(cf.~\cite{Audin.torus}).  The germs of such equivariant neighborhoods
are classified by their orbit types and their stabilizer
subgroups. The standard Hamiltonian action on $(\bR^4,\omega_0)$ has
point and torus fibers, and circle orbits whose stabilizer subgroups
are $G_{(1,0)}$ or $G_{(0,1)}$.  The stabilizer subgroups for point
and torus orbits are unique (equal to the torus and the identity,
respectively).  Meanwhile, an automorphism $A$ of the torus that is
acting on the $4$-manifold changes a stabilizer subgroup for a circle
orbit from $G_v$ to $G_{A^{T}v}$, allowing any stabilizer subgroup to
be achieved by an automorphism of the torus.
\end{proof}

\begin{obs} \label{O:stabilizersubgps} 
The moment map image of the neighborhood of a fixed point in a 
toric manifold is a convex sector 
bounded by rays with primitive integral tangent vectors $u,v$ such that 
the determinant $\lvert u v \rvert$ has norm $1$, and conversely all such
sectors appear as moment map images of $\bR^4$.
Meanwhile, the moment map image of the neighborhood of an orbit with stabilizer
subgroup $G_v$ is a neighborhood of a point in a half plane whose
boundary has $v$ as its tangent vector.
\end{obs}
Lemmas~\ref{L:stdtoricmodeluniqueness} and~\ref{L:AffineMoment} tell us that
Example~\ref{E:stdR4} and its variants induced by automorphisms of the
torus provide a complete set of local models for the neighborhood of
an orbit in a toric symplectic manifold, and that these models are
distinguished by their moment map images.

We now turn to the question of what toric structures can look like in the
neighborhood of a component of the vanishing locus in a
near-symplectic manifold. 

\begin{prop} \label{P:stdvanlocusuniq} Each component $C$ of the
  vanishing locus in a toric near-symplectic manifold
  $(X,\omega,\sigma)$ has an open torus-invariant neighborhood $N$ and
  an orbit-preserving map $\phi: (N,\omega,\sigma) \rightarrow (S^1
  \times \bR^3,\omega_1,\sigma_1)$ which is a smooth symplectic
  embedding of $N \setminus C$ that maps $C$ to $S^1 \times 0$ and is a
  homeomorphism onto its image.  Again, up to an automorphism of the
  torus, this embedding can be taken to be equivariant.
\end{prop}

\begin{proof}
Whenever we use Cartesian coordinates $(x,y,z)$ on $\bR^3$, then we
will freely also use cylindrical coordinates $(r,\theta,z)$ on
$\bR^3$ without further ado.

Following the discussion of smooth torus actions in
Section~\ref{S:LocalModels} (based on~\cite{OrlikRaymond}), we know
that, up to an automorphism of $T^2$, we can choose coordinates
$\alpha \in S^1$, $(x,y,z) \in \bR^3$ on a neighborhood $N$ of $C$
so that $C$ is $\{(\alpha,x,y,z) \mid x=y=z=0\}$ and so that the
action is $\sigma((t_1,t_2), (\alpha,r,\theta,z)) =
(\alpha+t_1,r,\theta+t_2,z)$.

Suppose $\mu = (\mu_1, \mu_2): N \rightarrow \bR^2$ is a moment map
for $\sigma$ on $N$.

The orbit space is 
homeomorphic to the upper half plane which, as 
in Example~\ref{E:stdvanishinglocus}, we parameterize 
by coordinates $X,Y$ so that the projection $\pi$ to the orbit space is
given by $X=z$ and $Y=r^2/2$.
(Note that we are using $(X,Y)$ as coordinates on this copy of $\bR^2$ to
distinguish them from $(x,y)$ which are coordinates on $N$.)

Then the moment map $\mu$ factors through the orbit space, $\mu=p\circ
\pi$ where $p = (p_1,p_2) : H \rightarrow \bR^2$. The map
$p$ is smooth and is an immersion on $H \setminus (0,0)$.
Because the isotropy subgroup of the circle orbits is $G_{(1,0)}$, the map
$p$ sends both the
positive and negative $X$-axes to straight lines of slope $0$. The
fact that $\mu$ is a moment map for the given action then means that
$\omega = dp_1 \wedge dq_1 + dp_2 \wedge dq_2$, where $q_1 = \alpha$
and $q_2 = \theta$.

Since we can freely translate the image of a moment map and can apply the
torus automorphism 
$\left(\begin{smallmatrix} -1 & 0 \\ 0 & 1 \end{smallmatrix}\right)$
without changing the isotropy subgroup, we assume
without loss of generality that $p$ maps the origin to the origin and
the positive $X$-axis to the positive $p_1$-axis. Consider a small
semicircular arc in $H$ starting on the positive $X$-axis and ending
on the negative $X$-axis, avoiding $(0,0)$. This arc is mapped by $p$
to a path starting on the positive $p_1$-axis and ending either on the
negative or positive $p_1$-axis, avoiding $(0,0)$, completing a total
of $k$ half-rotations, for some positive integer $k$ which is
independent of the choice of arc. We claim that $k=2$. In other words,
$p$ also maps the negative $X$-axis to the positive $p_1$-axis and is
injective on the interior of $H$.

We prove this claim by means of the following calculations.

We express $\omega$ in local coordinates near a point on $C$ as a
map $\omega: \bR^4 \rightarrow \bR^6$, where the $\bR^4$-coordinates
are $(\alpha,x,y,z)$ and the $\bR^6$ coordinates are the coefficients
of $(d\alpha \wedge dx, d\alpha \wedge dy, d\alpha \wedge dz, dx
\wedge dy, dx \wedge dz, dy \wedge dz)$. 
We compute differentials:
\begin{align*}
  dp_1 &= x \dponedY dx + y \dponedY dy + \dponedX dz \\
  dp_2 &= x \dptwodY dx + y \dptwodY dy + \dptwodX dz \\
  dq_1 &= d\alpha \\
  dq_2 &= \frac{1}{r^2} (-y dx + x dy) = - \frac{y}{2} \frac{1}{Y} dx
  + \frac{x}{2} \frac{1}{Y} dy.
\end{align*}
Then because
\begin{align*}
     dp_1 \wedge dq_1 + & dp_2 \wedge dq_2
       = - x \dponedY d\alpha \wedge dx - y \dponedY d\alpha \wedge dy -
       \dponedX d\alpha \wedge dz \\
       & + \dptwodY dx \wedge dy + \frac{y}{2} (\frac{1}{Y} \dptwodX) dx
       \wedge dz - \frac{x}{2} (\frac{1}{Y} \dptwodX) dy \wedge dz,
\end{align*}
we have the following expression for $\omega$ as a map to $\bR^6$:
\be
\omega(\alpha,x,y,z) = (- x \dponedY, - y \dponedY, -\dponedX, \dptwodY,
           \frac{y}{2} (\frac{1}{Y} \dptwodX),
           - \frac{x}{2} (\frac{1}{Y} \dptwodX))
\ee
The derivative $D \omega$ is simply the $6$-by-$4$ matrix of partial
derivatives of this function.  The claim will follow from the fact that,
on the vanishing locus $Z_\omega=\{x=y=z=0\}$ where $X=Y=0$,
this matrix must have rank $3$, together with the requirements
that $\omega$ must be everywhere well-defined and equal to $0$ 
on $Z_\omega$.
For $\omega$ to be well-defined, $\dptwodX$ must be divisible by
$Y$, i.e. 
\[ \dptwodX = Y f(X,Y)\]
for some smooth function $f$, and in particular 
\[ \dptwodX = 0 \ \ \ {\rm when}\ \ X=Y=0. \]
That $\omega$ must vanish on $Z_\omega$ implies
\[ \dponedX = \dptwodY=0   \ \ \ {\rm when}\ \ X=Y=0. \]
Now we compute $D \omega$, noting that
everything is $\alpha$-invariant, that for any function $h(X,Y)$ we have
\be 
\frac{\del h}{\del x} = x \frac{\del h}{\del Y}, \ \ 
\frac{\del h}{\del y} = y \frac{\del h}{\del Y}, \ \ 
\frac{\del h}{\del z} = \frac{\del h}{\del X},
\ee
and that
\be
\frac{\del^2 p_2}{\del X \del Y}  = \frac{\del}{\del Y} ( Y f(X,Y))
= f(X,Y) + Y \frac{\del f}{\del Y}
\ee
so that
\[ \frac{\del^2 p_2}{\del X \del Y}(0,0) = f(0,0). \]
Thus, on $Z_\omega$ where $X=Y=0$, we get:
\[ D \omega = \begin{bmatrix}
0 & - \dponedY & 0 & 0 \\
0 & 0 & - \dponedY & 0 \\
0 & 0 & 0 & -\frac{\del^2 p_1}{\del X^2} \\
0 & 0 & 0 & f \\
0 & 0 & \frac{1}{2} f & 0 \\
0 & -\frac{1}{2} f & 0 & 0
\end{bmatrix}
\]
Therefore, to have rank $3$, we need either that
\[ f(0,0) = \frac{\del^2 p_2}{\del X \del Y}(0,0) \neq 0 \]
or that
\[ \frac{\del^2 p_1}{\del X^2}(0,0) \neq 0 \quad \text{and} \quad
\dponedY(0,0) \neq 0.\]

We now complete the proof of the claim by contradiction. Suppose $k
\neq 2$. Note that $p_2$ is a real-valued function of two variables
mapping $(0,0)$ to $0$ with $(0,0)$ as a critical point. (Since
$p_2$ is smooth on $H$ we may extend its domain to an open
neighborhood of $H$.) If $k=1$ then $p_2(X,Y) > 0$ for all $Y > 0$,
while $p_2(X,0) = 0$ for all $X$, so $(0,0)$ must be a degenerate
critical point for $p_2$. Also if $k > 2$ then $(0,0)$ must be a
degenerate critical point, and so in either case the Hessian of
$p_2$ at $(0,0)$ must be singular. Since $\dptwodX = Y f(X,Y)$, we
know that $\frac{\del^2
  p_2}{\del X^2}(0,0) = 0$, so the Hessian being singular implies that
at $(0,0)$ we have $\frac{\del^2 p_2}{\del X \del Y} = 0$.

Now consider $p_1$. If $k=1$ then, along the $X$-axis, $p_1(X,0)$ is
an increasing function of $X$ whose first derivative vanishes at
$X=0$, so its second derivative must also vanish at $X=0$,
i.e. $\frac{\del^2 p_1}{\del X^2} = 0$ at $(0,0)$.  On the other hand
if $k> 2$, consider $p_1$ evaluated along semicircles centered at
$(0,0)$ and of radius $\epsilon$, as $\epsilon \ra 0$.  Define
$\Gamma_\epsilon:=\{(X,Y) \,|\, X^2+Y^2=\epsilon^2 \ {\rm and}
\frac{\del^2 p_1}{\del X^2}=0\}$.  On each semicircle $\frac{\del^2
p_1}{\del X^2}$ takes on both positive and negative values, so for
each $\epsilon>0$, $\Gamma_\epsilon$ is nonempty.  Since $\frac{\del^2
p_1}{\del X^2}=0$ is a closed condition, $\cup_\epsilon
\Gamma_\epsilon$ is a closed subset of $H$, implying that $\Gamma_0$
is nonempty.  Therefore, $\frac{\del^2 p_1}{\del X^2}=0$ at $(0,0)$.

Thus we must have $k=2$, which establishes the claim.

The above claim means that these coordinates $(p_1,q_1,p_2,q_2)$ behave,
topologically, exactly the same as the $(p_1,q_1,p_2,q_2)$ coordinates
in Example~\ref{E:stdvanishinglocus}. Thus, identifying these
coordinates here with the corresponding coordinates in
Example~\ref{E:stdvanishinglocus} gives the desired homeomorphism which can
fail to be smooth only on the vanishing locus $p_1=p_2=0$.
\end{proof}

\section{Toric near-symplectic manifolds}
\label{S:ToricAndIAS}

Recall that a toric near-symplectic manifold is a near-symplectic manifold
equipped with an effective smooth torus action that is Hamiltonian on the 
symplectic locus.  As such, there are two relevant classifications that are
well-understood: that of $T^2$-manifolds in the smooth category and that
of toric manifolds in the symplectic category.

Orlik and Raymond classify
$T^2$-manifolds in terms of their weighted orbit spaces.
They first note that the orbit space is a surface with boundary
such that each point in the interior of the surface is the image of a torus
and each point on the boundary is the image of a lower-dimensional
fiber (circle or point).  Furthermore, the boundary is comprised of a
union of {\em edges}, the interiors of which parameterize circle orbits and
the endpoints of which are the images of point orbits.  The 
{\em weighted orbit space} of a $T^2$-manifold is then the oriented orbit 
space together
with a labeling of the edges by the corresponding stabilizer subgroups.

Meanwhile, Delzant's Theorem classifies 
closed toric manifolds in terms of their moment map 
images~\cite{Delzant.moment}, which are polygons (and in fact are weighted
orbit spaces because the tangent vectors to their edges encode the
stabilizer subgroups).

We cannot apply Delzant's classification to the symplectic locus of a
toric near-symplectic manifold because it is noncompact, thereby
allowing the preimage of a point in the moment map image to be
disconnected.  To accommodate this feature we introduce integral
affine surfaces, which are essentially weighted orbit spaces with an
induced geometry that encodes the essential structure of the moment
map.  We then determine the extent to which one can classify such
manifolds in terms of their orbit spaces equipped with integral affine
structures (Theorem~\ref{T:ToricUniqueness}).  For a more leisurely
discussion of this matter in the symplectic case the reader can
consult~\cite{Symington.4from2}.

Consider a toric near-symplectic manifold $(X,\omega,\sigma)$ and the
projection to its orbit space, $\pi:X\ra B$.  Let $F$ be the discrete
subset of $\bdy B$ that is the image of the vanishing locus.  Then the
restriction of the moment map $\mu$ to the symplectic locus factors
through an orientation preserving immersion $\Phi:B\setminus F\ra
\bR^2$ which extends to a continuous map $\overline{\Phi}$ such that
$\mu=\overline{\Phi}\circ\pi$.

The immersion $\Phi$ 
induces a geometric structure on the orbit space whose local isometries
are elements of ${\rm Aff}(2,\bZ)$.  Henceforth, unless otherwise noted,
we assume that \lq\lq surface\rq\rq\ means \lq\lq surface with possibly
nonempty boundary.\rq\rq

\begin{defn} \label{D:IAS}
  An {\em integral affine structure} on a surface is
  a maximal atlas of charts to a sector of $\bR^2$ whose boundary
  rays have rational slope, with transition functions in
  ${\rm Aff}(2,\bZ)$. The {\em standard integral affine structure} on
  $\bR^2$, denoted $\cA_0$, is the atlas containing the identity map.
  Two {\em integral affine surfaces} $(B,\cA)$ and
  $(B',\cA')$ are {\em isomorphic} if there
  exists a homeomorphism $\phi:B\ra B'$ such that $\phi^*\cA'=\cA$.
\end{defn}
Note that the homeomorphism $\phi$ is a diffeomorphism on the complement
of the vertices.

\begin{defn}
  A {\em vertex} on the boundary of an integral affine surface is a point
  that maps, via a chart, to a vertex of a sector in $\bR^2$, while an
  {\em edge} is the closure of a connected component of the boundary
  minus its vertices.
  The boundary of an integral affine surface is {\em right polygonal}
  if every vertex has a neighborhood that is isomorphic to a 
  neighborhood of the origin in the quadrant $(Q,\cA_0)\subset(\bR^2,\cA_0)$.
\end{defn}

For near-symplectic manifolds equipped with a toric structure, the local
model near a component of the vanishing locus 
(Example~\ref{E:stdvanishinglocus}) dictates
the geometric structure in the base near the image of such a component
(Figure~\ref{F:onefold}).
This inspires the following:
\begin{defn} \label{D:IASEdgeFolds} The boundary of an integral affine
  surface has {\em edge folds} if there is a set of points $F \subset \bdy B$,
  called {\em fold points}, such that each point $p \in F$ has a neighborhood 
  $U$ and a homeomorphism $\phi: U
  \rightarrow V$ to a neighborhood $V$ of $(0,0)$ in the upper
  half-plane $H = \{(X,Y) \mid Y \geq 0\}$, such that the integral
  affine structure on $U \setminus p$ is $\phi^* \psi^*
  \mathcal{A}_0$, where $\psi: (X,Y) \mapsto (X^2 - Y, 2XY)$ and
  $\mathcal{A}_0$ is the standard integral affine structure on
  $\bR^2$.  

  If the set of fold points is nonempty, we always indicate these singularities
  explicity, so the integral affine surface whose boundary has edge folds
  is a triple $(B,\cA,F)$.
  Two surfaces with such structures,
  $(B,\mathcal{A},F)$ and $(B',\mathcal{A}',F')$, are equivalent if there 
  exists a
  homeomorphism from $B$ to $B'$, taking $F$ bijectively to
  $F'$, which is an integral affine equivalence of $(B\setminus
  F,\mathcal{A})$ and $(B' \setminus F', \mathcal{A}')$.
\end{defn}

The orbit space of a toric near-symplectic manifold acquires the 
structure of an integral
affine surface whose boundary is right polygonal with edge folds
by pulling back
to $B\setminus F$ the standard integral affine structure via the immersion 
$\Phi$.  
Note that Lemma~\ref{L:AffineMoment} implies 
that if two toric near-symplectic manifolds are equipped with
actions that differ by an automorphism of the torus, then the integral
affine structures on their orbit spaces are isomorphic.

\begin{defn} \label{D:intaffbase}
Given a toric near-symplectic manifold $(X,\omega,\sigma)$ with
moment map $\mu$ and orbit space projection $\pi$, there is a unique
map $\overline\Phi$ defined by $\mu=\overline\Phi\circ\pi$.  The
local models for the moment map near any orbit imply that
$\overline\Phi$ is an immersion on the complement of a set of
isolated points which we denote by $F$.  Let
$\Phi=\overline\Phi|_{B\setminus F}$, where $B=\pi(X)$.  The {\em
integral affine base} of the fibration defined by $\pi$ and $\mu$ is
then $(B,\Phi^*\cA_0,F)$, an integral affine surface whose boundary
is right polygonal with edge folds.
\end{defn}

\begin{thm} \label{T:ToricBases}
  Suppose $(B,\cA,F)$ is an
  integral affine surface whose boundary is right polygonal boundary 
  with edge folds.  
  Then $(B,\cA,F)$ is the integral affine base of a toric near-symplectic
  manifold if and only if there is
  an integral affine immersion $\Phi:(B\setminus F,\cA)\ra(\bR^2,\cA_0)$
  that extends to a continuous map $\overline{\Phi}:B\ra\bR^2$
  so that each point in $F$ has a neighborhood $N$ on which 
  $\overline\Phi|_{N\setminus\bdy B}$ is injective and
  $\overline\Phi|_{(N\setminus F)\cap\bdy B}$ is two-to-one and linear.
\end{thm}
Notice that there is no requirement here that $B$ be compact.  

\begin{proof}
The \lq\lq only if\rq\rq\ direction follows directly from the
factorization of the moment map for a toric near-symplectic manifold
mentioned just before Definition~\ref{D:IAS} and the local models in
Examples~\ref{E:stdR4} and~\ref{E:stdvanishinglocus}, together with
Lemma~\ref{L:stdtoricmodeluniqueness} and
Proposition~\ref{P:stdvanlocusuniq}, invoking automorphisms of the
torus as needed.  The \lq\lq if\rq\rq\ direction follows by
construction.

Fix a torus $T^2$ with cyclic coordinates $q=(q_1,q_2)$ of period $2\pi$
and consider the manifold $B\times T^2$ with the smooth torus
action $t\cdot(x,q)=(x,q+t)$ where $t=(t_1,t_2)$.
Equip $B\setminus F$ 
with local coordinates $(p_1,p_2)=\overline{\Phi}(x)$,  $x\in B$.
Then $dp\wedge dq$ is a symplectic form on $(B\setminus F)\times T^2$ with
respect to which the action is Hamiltonian.  Notice that the moment map for
this action amounts to forgetting the cyclic coordinates.

Take the quotient of $B\times T^2$ by identifying points as follows:
if $x$ belongs to a vertex in the boundary of $(B,\cA,F)$ then identify
$(x,q)$ and $(x,q')$ for all $q,q'$; if $x$ belongs to the interior of
an edge $E$ of $B$, then identify $(x,q)$ and $(x,q')$ whenever 
$t\cdot(x,q)=(x,q')$ for some $t\in G_v$ where $v$ is a primitive integral
tangent vector to $\overline\Phi(E)$.  
Example~\ref{E:stdR4}, together with automorphisms of the torus, and the fact 
that $(B\setminus F,\cA)$ is an integral affine surface with
right polygonal boundary, assure that the quotient is a smooth manifold with
a smooth torus action, and that the symplectic structure (on the
symplectic locus) descends to 
the quotient, making the action Hamiltonian there.
Call this manifold $X$ and let $\omega$ be the symplectic structure defined
on the symplectic locus.  Let $\pi:X\ra B$ be projection induced on $X$
by the projection $B\times T^2\ra B$ that merely forgets the $T^2$ factor.

Now, for any point $p$ in $F$, let $N$ be a small neighborhood of $b$ 
satisfying the hypotheses.  Let $N'=N\setminus p$.  Then by construction, 
up to an automorphism
of the torus, $(\pi^{-1}(N'), \omega|_{\pi^{-1}(N')})$ 
embeds symplectically and equivariantly
into the complement of the vanishing locus in 
Example~\ref{E:stdvanishinglocus}.  Let $\phi$ be the continuous extension
of this embedding to $\pi^{-1}(N)$.  Pulling back via $\phi$ 
the near-symplectic structure
$\omega_1$ and the torus action of Example~\ref{E:stdvanishinglocus} we get a
toric structure on $\pi^{-1}(N)$ that is compatible with the toric structure
on the symplectic locus of $X$.  Doing this for each $p\in F$ we complete the
construction.
\end{proof}

The question of what integral affine surfaces immerse in $(\bR^2,\cA_0)$ is
being investigated by the second author.  However, there is one case in
which one can be assured of the existence of an integral affine immersion.
\begin{lemma}
\label{L:SimplyConnected}
If $(B,\cA)$ is an integral affine surface such that $B$ has trivial 
fundamental group, then there exists an integral affine immersion
$(B,\cA)\ra(\bR^2,\cA_0)$.
\end{lemma}
\begin{proof}
$(B,\cA)$ is a manifold locally modeled on $(\bR^2,\cA_0)$.  As such,
its developing map (cf.~\cite{Thurston.3manifolds}) 
$D:(\widetilde B,\Psi^*\cA)\ra(\bR^2,\cA_0)$ from the 
universal cover (whose covering map is 
$\Psi:(\widetilde B, \Psi^*\cA)\ra(B,\cA)$) is an integral affine immersion.
But since $B$ is contractible, $\widetilde B=B$ and $\Psi=Id$. 
\end{proof}

We can now state and prove the generalization of Delzant's theorem for
toric near-symplectic manifolds.

\begin{thm} \label{T:ToricUniqueness} An integral
  affine surface $(B,\cA,F)$, with possibly nonempty right polygonal boundary
  with edge folds, together with an integral affine immersion 
  $\Phi:(B\setminus F,\cA)\ra(\bR^2,\cA_0)$, 
  determines a toric near-symplectic 
  manifold uniquely up to equivariant homeomorphisms that, restricted to the
  symplectic locus, are symplectomorphisms.  The integral affine base
  $(B,\cA,F)$ by itself determines the toric near-symplectic
  manifold up to an orbit-preserving homeomorphism that is a 
  symplectomorphism on the symplectic locus.
\end{thm}

\begin{proof}
Without loss of generality, we restrict to open covers such that each component
of the vanishing locus belongs to a unique open set -- so the gluing of
open sets occurs only on the symplectic locus.
Use Theorem~\ref{T:ToricBases} to 
construct a toric near-symplectic manifold $(X,\omega,\sigma)$ 
whose integral affine base is 
$(B,\cA,F)$ with projection $\pi:X\ra B$.
Let $\cS$ be the sheaf on $B$ that associates
to any open set $U\subset B$ the group of equivariant 
symplectomorphisms of $(\pi^{-1}(U),\omega_{\pi^{-1}(U)})$.  
The sheaf cohomology group $H^1(B,\cS)$
encodes the transition maps between such torus invariant neighborhoods and
classifies, up to isomorphism, toric manifolds whose moment map $\mu$ 
factors as $\Phi\circ\pi$.  

Our goal is to show that $H^1(B,\cS)$ is trivial.  
To do so we follow the argument put forth 
in~\cite{LermanTolman.orbifold} (Prop. 7.3), which the reader may consult
for further details.
Since the equivariant symplectomorphisms are all time one flows of Hamiltonian
vector fields, $\cS$ fits into an exact sequence
of sheaves,
\be
\label{E:shortexactseq}
0\ra\cL\ra\cC^\infty\ra\cS\ra 0,
\ee
where $\cC^\infty$ is the sheaf of smooth torus invariant functions on 
$X$ (which is equivalent to the sheaf of smooth functions on $B$),
\footnote{In~\cite{LermanTolman.orbifold} 
the sheaf $C^\infty$ is isomorphic to
the sheaf of continuous functions on $B$ that lift to smooth functions on 
the total space 
because the authors are working with orbifolds.}
and $\cL$ is the sheaf of torus invariant
functions whose time one flows generate the identity map. 
This short exact sequence induces the long exact sequence
\be
\label{E:longexactseq}
\cdots \ra H^i(B,\cL) \ra H^i(B,\cC^\infty) \ra H^i(B,\cS) \ra
 H^{i+1}(B,\cL) \ra \cdots,
\ee
so it suffices to show that $H^1(B,\cC^\infty)=H^2(B,\cL)=0$

The sheaf $\cL$ is isomorphic to the sheaf of locally constant functions
on $B$ with values in $\bR\times \Lambda^*$ where $\Lambda^*\cong\bZ^2$ 
is the lattice in 
the Lie algebra $\mathbf{t}^*$ consisting of covectors $\xi$ whose 
infinitesimal
vector field has the identity map as its time one flow.  Indeed, the elements
of $\cL$ are precisely the functions $f_{(c,\,\xi)}$ defined by 
$f_{(c,\,\xi)}(b)= c+\left<\xi,\Phi(b)\right>$ 
for each $(c,\xi)\in\bR\times\Lambda^*$.
Because $B$ is locally contractible, we have isomorphisms 
between sheaf cohomology with coefficients in $\bR$ and $\bZ$, and both 
deRham and singular cohomology.  Consequently,  $H^2(B,\cL)=0$ because
\be
\label{E:Linear}
H^i(B,\cL)\cong H^i(B,\bR)\times H^i(B,\bZ)\times H^i(B,\bZ),
\ee
and the fact that $B$ immerses in $\bR^2$ implies 
$H^2(B,\bZ)\cong H^2(B,\bR)=0$.
Meanwhile, because $\bC^\infty$ is a fine sheaf, $H^i(B,\cC^\infty)=0$
for all $i>0$.  

A different choice of integral affine immersion would yield a new toric 
manifold with moment map
$\Phi'=\Psi\circ\Phi$ for some $\Psi\in{\rm Aff}(2,\bZ)$.
If $A$ were the linear part of $\Psi$, then by Lemma~\ref{L:AffineMoment} 
and the uniqueness proved above, the corresponding toric manifold 
would be $(X,\omega,\sigma')$ with 
$\sigma'=\sigma\circ(A^{T}\times Id)$, which is indeed orbit preserving 
symplectomorphic to $(X,\omega,\sigma)$.
\end{proof}

The fact that our maps may not be smooth across the vanishing loci is
a reflection of the fact~\cite{Honda.local} that Moser's method for
near-symplectic forms near their vanishing loci does not give
smoothness at the vanishing loci.  Thus it is really a feature of the
germ of the toric structure on a component of the vanishing locus, rather
than a feature of how we glue in a neighborhood of the component.

\section{Making $T^2$-manifolds near-symplectic}
\label{S:GlobalTorusAction}

\subsection{Simply connected $T^2$-manifolds}

Equipped with Orlik and Raymond's technology~\cite{OrlikRaymond} for
understanding torus actions in the smooth category, as outlined in
Section~\ref{S:ToricAndIAS}, we now investigate the problem of finding
near-symplectic structures adapted to given smooth torus actions.
The goal is to prove Theorems~\ref{T:SimplestThm}, \ref{T:HamiltActions}
and Proposition~\ref{P:SimplyConnectedT2}.
In light of Theorem~\ref{T:ToricUniqueness}, proving 
Theorems~\ref{T:SimplestThm} and \ref{T:HamiltActions}
amounts to realizing given weighted orbit spaces as 
integral affine surfaces with appropriate boundaries 
isometrically immersed in $(\bR^2,\cA_0)$.

We separate 
Theorem~\ref{T:SimplestThm} into three propositions (\ref{P:SimplestThmI},
\ref{P:SimplestThmII}, and~\ref{P:SimplestThmIII}), 
the first of
which is:
\begin{prop}
\label{P:SimplestThmI}
  Every locally toric fibration of a closed simply connected near-symplectic
  manifold is toric.
\end{prop}

\begin{proof}
Let $\pi:(X,\omega)\ra(B,\cA,F)$ be the locally toric fibration.
By Theorem~\ref{T:ToricBases}, $(X,\omega)$ is toric if there is
an integral affine immersion $\Phi:(B\setminus F,\cA)\ra(\bR^2,\cA_0)$.

Because any loop in $B$ lifts to a loop in $X$, $X$ being simply connected
implies $B$ is simply connected.  Because a sphere cannot admit an integral 
affine structure, $B$ must be homeomorphic to a disk.
Finally, by Lemma~\ref{L:SimplyConnected} we know that $(B\setminus F,\cA)$
does immerse isometrically in $(\bR^2,\cA_0)$. 
\end{proof}

Before proving the rest of Theorem~\ref{T:SimplestThm}
we introduce language to describe the image of the boundary
of an immersed integral affine surface.

\begin{defn}
\label{D:SlopeLists}
A list of slopes $s_1,\ldots, s_k$, or {\em slope list}, is 
{\em right polygonal} 
if $s_j\in\bQ\cup\infty$ for each $j$ and, writing $s_j=\frac{m_j}{n_j}$
as a reduced fraction (with $m_j=1,n_j=0$ if $s_j=\infty$),  
if the determinant
$\bigl(\begin{smallmatrix} m_j & m_{j+1} \\ n_j 
& n_{j+1} \end{smallmatrix} \bigr)$ has norm $1$
for each $j\le n-1$.
When the indices of a list of $k$ slopes are understood mod $k$, the slope 
list is {\em cyclic}.
\end{defn}

\begin{defn}
\label{D:PolygonalPaths}
A {\em polygonal path with folds} is 
a piecewise linear map $\gamma:[0,n]\ra\bR^2$, $n\in\bN$, well defined up to
reparametrization relative endpoints on the subintervals $[j-1,j]$, 
$j=1,\ldots N$, such that 
\ben
\item the image $e_j$ of $\gamma|_{[j-1,j]}$, $j\in\{1,\ldots N\}$, 
has constant slope (and is called an {\em edge}) and
\item the slopes of $e_j$ and $e_{j+1}$ are different for each 
$j=1,\ldots, N-1$.
\een
The {\em vertices} are the images of the integral points, $\gamma(j)$, while
the {\em fold points} are the interior points of edges at which $\gamma$
is not smooth, i.e. where the oriented tangent vector changes direction.

A polygonal path (possibly with folds) 
is {\em right polygonal} if each vertex $\gamma(j)$ has a
neighborhood that contains no fold points and in which 
oriented tangent vectors $v_j$ and
$v_{j+1}$ to the edges $e_j,e_{j+1}$ 
have determinant $\lvert v_j v_{j+1} \rvert = 1$.
\end{defn}
The definitions are analogous for {\em polygonal loops} with the usual
additional stipulation that $\gamma(0)=\gamma(n)$.
Notice that polygonal paths and loops are oriented in accordance with
the standard orientation of the domain interval $[0,n]$ -- and the indexing
of the edges and vertices.

Given a (cyclic) slope list, we say that a polygonal path
(or loop) with folds 
{\em realizes} this list if the list of slopes of the edges, ordered according
to the orientation,  equals the given slope list.

\begin{defn} \label{D:Embedded}
  Abusing correct terminology considerably, let us say that a polygonal
  path with folds is {\em embedded} if the only intersections between distinct
  edges occur when two consecutive edges meet at a vertex and if there are no
  triple points.
\end{defn}
In other words, an edge with a fold obviously intersects itself, but
this is the only type of self-intersection allowed and such an edge
may not intersect itself too much. In particular a given edge may have
at most two folds. To justify the use of the term \lq\lq embedded\rq\rq, note
that if we remove all double points and take the closure of the
remainder, we get an honestly embedded polygonal path.

\begin{defn} \label{D:PositiveTurning}
  Given two slopes $s_1, s_2 \in \bQ \cup \{\infty\}$, define the {\em
  positive angle} from $s_1$ to $s_2$ to be the angle
  $\alpha(s_1,s_2) \in [0,\pi)$ from a line of slope $s_1$ to a line
  of slope $s_2$ measured counterclockwise.  Given a list of slopes
  $s_1, \ldots, s_n \in \bQ \cup \{\infty\}$, define the {\em positive
  turning} for the list to be $T := \alpha(s_1,s_2) + \ldots +
  \alpha(s_{n-1},s_n)$  
  or, if the list is cyclic, $T := \alpha(s_1,s_2) + \ldots +
  \alpha(s_{n-1},s_n) + \alpha(s_n,s_1)$.  Given a $T^2$-manifold
  with weighted orbit space $B$,
  the {\em positive turning} for a component of $\partial B$ 
  with weights $\{\frac{m_i}{n_i}\}$, $i=1,\ldots n$, is the 
  positive turning for the cyclic
  slope list $\{\frac{m_i}{n_i}\}$, $i=1,\ldots n$.
\end{defn}
Note that, for any list, $T=k\pi$ for some nonnegative integer $k$ and
that $T=0$ if and only if the boundary component has only one edge.

\begin{defn} \label{D:TotalTurning}
  Given a polygonal path with folds, define the {\em total turning} to
  be the sum of the amounts of counterclockwise turning of the tangent 
  vectors at the vertices (between $0$ and
  $\pi$) and at the folds (always $-\pi$). If the path is closed, we include
  the turning at the initial vertex (which is also the final vertex).
\end{defn}

The following is an immediate consequence of
Definitions~\ref{D:Embedded}, \ref{D:PositiveTurning} 
and ~\ref{D:TotalTurning}.

\begin{lemma} \label{L:PosTurningMax} Given a slope list with positive
  turning $T$ and a right polygonal loop with folds realizing the list
  with total turning $T'$, then  $T' = T - \abs{F}\pi$,
  where $\abs{F}$ is the number of folds. An embedded right polygonal loop
  with folds in $\bR^2$ always has total turning $\pm 2\pi$.
\end{lemma}

\begin{lemma}
\label{L:PolygonalPath}
Given a right polygonal slope list $s_1,\ldots s_k$, $k\ge2$, with $s_1=0$ 
and $s_2=\infty$, and a point $(x_0,y_0)\in\bR^2$ with $x_0,y_0>0$, 
there exists an embedded 
right polygonal path with folds that realizes the slope list, 
has initial endpoint $(0,0)$ and final endpoint $(x_0,y_0)$, and whose
first edge $e_1$ has oriented tangent vector $(1,0)$.
\end{lemma}

\begin{proof}
Let $a,b,\delta$ and $\epsilon$ be indeterminates that will be fixed later.
Let $e_1 = [0,a] \times \{0\}$ and let $e_2$ be an
edge that starts with $\{a\} \times [0,b]$, has a fold at $(a,b)$ and
then doubles back a distance $\delta$, ending at $(a, b-\delta)$. Now
turn counterclockwise onto an edge $e_3$ of slope $s_3$ and length
$\epsilon$. Henceforth, turn counterclockwise from $e_i$ onto an edge
$e_{i+1}$ of slope $s_{i+1}$ and length $\epsilon$ if this does not
require moving in the negative $x$ direction. Otherwise, put a fold at
the current end of $e_i$ and extend $e_i$ by doubling back a distance
$\delta$, then turn counterclockwise onto $e_{i+1}$, of slope
$s_{i+1}$ and length $\epsilon$. Continue up to $e_k$.  Given any
preassigned $\lambda > 0$, there exists a choice of $\epsilon$ and
$\delta$ ($\epsilon$ small and $\delta$ much smaller) such that, for
any $a, b > 0$, the entire path is embedded and ends at $(a+a_0,
b+b_0)$, where $0 \leq a_0 < \lambda$ and $-\lambda < b_0 <
\lambda$. Then, if we choose $a = x_0 - a_0$ and $b = y_0 - b_0$, we
can arrange that the path ends at $(x_0,y_0)$.
\end{proof}

Notice that a slope list $s_1,\ldots, s_n$ can be realized by a 
right polygonal path if and only if one can also
realize the slope list
$s_1',\ldots, s_n'$ where the reduced fractions $m_i/n_i$
and $m'_i/n'_i$ representing the slopes satisfy
$(n'_i,m'_i)^T = A (n_i,m_i)^T$
for some fixed $A\in GL(2,\bZ)$.

\begin{lemma} \label{L:SimpleLoop}
Any cyclic reduced slope list with positive turning $T \geq 2\pi$ can
be realized by an embedded right polygonal loop with folds.
\end{lemma}

\begin{proof}
Let $s_1, \ldots, s_n$ be the given list. Without loss of generality
assume that $s_1=0$ and, for each $i$, $s_i \neq \infty$. 
Let $p$ be the smallest integer between $2$ and $n$ such that
$s_p \neq 0$ and such that the cyclic slope list $s_1,s_p,\ldots, s_n$ has
total positive turning $2\pi$.  (Notice that this new slope list may
fail to be right polygonal at the intersection of the $p^{\rm th}$
and $1^{\rm st}$ edges.) Without loss of generality we can
also assume that $s_p > 0$.
If $p=2$ then it is standard to
construct the desired loop, a closed convex polygon without folds. Now
assume $p > 2$. 

Figure~\ref{F:SimpleLoop} illustrates the following
construction. 
\begin{figure}
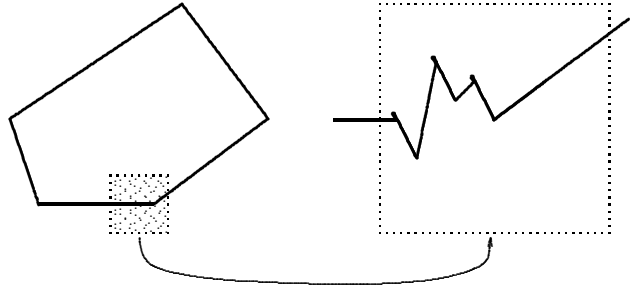

\begin{center}
\include{SimpleLoop}
\caption{An example construction of a right polygonal loop with folds.}
\label{F:SimpleLoop}
\end{center}
\end{figure}

First construct a closed convex polygon (without folds) with edges
$e_1,e_p \ldots, e_n$ realizing the cyclic slope list 
$s_1, s_p, \ldots, s_n$, with the
vertex joining $e_1$ to $e_p$ located at $(0,0)$.
Next, choose some $(x_0,y_0)$ on the interior of $e_p$ and 
replace the portion of $e_p$
connecting $(0,0)$ and $(x_0,y_0)$ by a right polygonal path
representing the slope list $s_1,\ldots, s_{p-1}$,
constructed as in Lemma~\ref{L:PolygonalPath}.
Note that in doing so the edge $e_p$ is shortened, its initial point 
becoming $(x_0,y_0)$, while $e_1$ is lengthened, its
final point becoming $(a,0)$ rather than $(0,0)$.
We can assure that this procedure yields an embedded polygonal loop
by choosing the point $(x_0,y_0)$ sufficiently close to $(0,0)$ and
the parameter $\lambda$ of Lemma~\ref{L:PolygonalPath} sufficiently small.
\end{proof}

The next part of Theorem~\ref{T:SimplestThm} that we prove is:
\begin{prop}
\label{P:SimplestThmII}
  Every closed simply connected
  $T^2$-manifold $(X,\sigma)$ with $b_2^+(X)>0$ is toric with respect to
  some near-symplectic structure.  
\end{prop}

\begin{proof}
The $T^2$ structure yields a slope list with positive turning $T$ for some
$T\ge 0$.  Assume for the moment that $T\ge 2\pi$.

Appealing to Lemma~\ref{L:SimpleLoop} and its proof, construct an
embedded right polygonal loop with folds representing the given slope
list, thereby defining a near-symplectic structure with respect to
which $(X,\sigma)$ is toric. This loop is then the right polygonal
boundary with edge folds of an integral affine disk.

If $T\le\pi$ any attempt to draw such an embedded polygonal loop
fails: all right polygonal paths representing the slope list fail to
close up.  It remains to show that $T\le\pi$ cannot occur.  First note
the $T=0$ cannot occur for any $T^2$-manifold because a slope list
with just one slope would define a weighted orbit space with just one
edge and one vertex.  With only one isotropy subgroup, it is
impossible to satisfy the condition that $\lvert u v \rvert = \pm 1$
were $G_u,G_v$ are the isotropy subgroups of circle orbits in the
neighborhood of the fixed point.

Now suppose $T=\pi$.  
Consider the weighted orbit space for $(X,\sigma)$ with edges 
$e_1,\ldots e_k$ having slopes $s_i=\frac{m_i}{n_i}$ for each
$i=1,\ldots k$, where the signs of $m_i,n_i$ are chosen so that 
$\lvert\begin{smallmatrix} n_i & n_{i+1} \\ 
m_i & m_{i+1}\end{smallmatrix}\rvert=1$ for each $i=1,\ldots k-1$.
(The fact that $T=\pi$ corresponds to that fact that, with these choices
of signs, $\lvert\begin{smallmatrix} n_k & n_{1} \\ 
m_k & m_{1}\end{smallmatrix}\rvert=-1$.)

Consider the $T^2$-manifold 
defined by the weighted orbit space that is a disk with weights
given by the cyclic slope list $s_1,\ldots s_k, s_{k+1}, s_{k+2}$ 
where $s_{k+1}=s_1$ and $s_{k+2}=\frac{ m_{k+2}}{n_{k+2}}$
with $m_{k+2},n_{k+2}$ satisfying
$\left(\begin{smallmatrix} n_{k+2} & n_{1} \\ 
m_{k+2} & m_{1}\end{smallmatrix}\right) =1$.
Call this new $T^2$-manifold $(X',\sigma')$.  Then, topologically,
$X'$ is obtained from $X$ by removing a copy of $S^1\times D^3$
and gluing in a $D^2\times S^2$.  Therefore $b_2^+(X')=b_2^+(X)+1$.

Because 
\be 
\left| \begin{array}{cc} 
n_k & -n_{k+1} \\ 
m_k & -m_{k+1} \end{array}\right| 
= 
\left| \begin{array}{cc} 
-n_{k+1} & n_{k+2}\\ 
-m_{k+1} & m_{k+2} \end{array}\right|
=
\left| \begin{array}{cc} 
n_{k+2} & n_{1} \\ 
m_{k+2} & m_{1} \end{array}\right|
=1,
\ee
the total turning of the new cyclic slope list is $2\pi$ and 
$X'$ admits a symplectic structure with respect to which $\sigma'$ is 
Hamiltonian.  Therefore $b_2^+(X')=1$, implying $b_2^+(X)=0$.
\end{proof}

Finally, the third part of Theorem~\ref{T:SimplestThm} is:
\begin{prop}
\label{P:SimplestThmIII}
  If $(X,\omega,\sigma)$ is a closed toric near-symplectic manifold that 
  is simply connected, then the vanishing locus $Z_\omega$ must have exactly  
  $\abs{Z_\omega}=b_2^+(X)-1$ components.
\end{prop}
Notice that Proposition~\ref{P:SimplestThmIII} is trivially true
when $\abs{Z_\omega}=0$ because the only toric symplectic manifolds are
$S^2\times S^2$ and $\bC P^2\#\overline{\bC P}^2$, all of which have $b_2^+=1$.

\begin{example}
\label{E:ConnectedSum}
Consider the moment map image of $\bC P^2\#\bC P^2$ shown in
Figure~\ref{F:FoldExamples}, for which one can easily verify
Proposition~\ref{P:SimplestThmIII}.  This toric near-symplectic
manifold is an equivariant connected sum, and the decomposition can be
performed via a symplectic cut~\cite{Lerman.SymplecticCuts} along the
preimage of the vertical line segment that connects the fold point and
the horizontal edge (cutting $X$ along that $3$-sphere and then
collapsing the circles on the resulting boundaries that are in the
kernel of the symplectic form).  Note that the preimage of the
vertical segment is indeed a $3$-sphere because the circles orbits
mapping to the endpoints of the segment have isotropy subgroups
$G_{(1,0)}$ and $G_{(0,1)}$, and $\lvert \begin{smallmatrix} 1 & 0 \\ 0
& 1 \end{smallmatrix}\rvert$ has norm $1$.
\end{example}

The essence of our proof of Proposition~\ref{P:SimplestThmIII}
is to decompose $(X,\omega)$ via symplectic cutting
as a connected sum of manifolds, 
each of which has no vanishing locus, and
show that the quantity $b_2^+(X)-\abs{Z_\omega}$ 
remains constant for all the manifolds involved.
However, in general a toric manifold may have a $3$-sphere on which a
one-dimensional subtorus acts freely, with one orbit being a component
of the vanishing locus i.e. a $3$-sphere whose moment map
image is a line segment with one endpoint being a fold point.
For instance, consider 
the $T^2$-manifold with weighted orbit space whose slope list is
$\infty, 0, -1, -2$. 
This manifold is diffeomorphic to $\bC P^2\#\bC P^2$, but is 
not an equivariant connected sum.
(We leave it to the reader to use the techniques of this section to draw
the moment map image for this action.) Here we thank Brett Parker for
asking us a question which brought the possibility of such examples to
our attention.

We circumvent this by recognizing that 
\ben
\item in light of Lemma~\ref{L:PosTurningMax}, 
$b_2^+(X)-\abs{Z_\omega}$ depends only on the underlying $T^2$-manifold and the
existence of a near-symplectic structure that makes the $T^2$-manifold
toric, but does not depend on the particular near-symplectic structure; and
\item  the quantity $b_2^+(X)-\abs{Z_\omega}$ is 
invariant under equivariant blowups.
\een

\begin{proof}[Proof of Proposition~\ref{P:SimplestThmIII}]
In this proof, let manifolds have more than one connected component.
Let $(X_i,\omega_i,\sigma_i)$, $i=1,\ldots,N$ be a sequence of 
toric near-symplectic manifolds such that 
$(X_1,\omega_1,\sigma_1)$ and $(X,\omega,\sigma)$ are equivalent as
$T^2$-manifolds, 
$\abs{Z_{i+1}}=\abs{Z_i}-1$, $\abs{Z_n}=0$, and
$(X_{i+1},\omega_{i+1},\sigma_{i+1})$ is obtained from
$(X_{i},\omega_{i},\sigma_{i})$ by 
\ben
\item choosing a convenient near-symplectic structure on the
$T^2$-manifold $(X_i,\sigma_i)$,
\item performing equivariant blow-ups of some connected 
component of $(X_i,\sigma_i)$ as
necessary, and then
\item performing a symplectic cut of that component 
along a $3$-sphere that contains a component of the vanishing locus.
\een

The proposition will be proved if we show that $b_2^+(X_i)-\abs{Z_i}=c_i$ 
where $c_i$ is the number of connected components of $X_i$. 
This is true for $X_N$ because $\abs{Z_N}=0$ and each component has
$b_2^+=1$.
All we need to do is show that we can always find a near-symplectic
structure such that the toric manifold with underlying $T^2$-manifold
$(X_i,\sigma_i)$, or a blow-up of it, 
can be equivariantly decomposed so as to reduce by 
$1$ the number of fold points in the boundary of its moment map image.
 
Given the toric manifold $(X_i,\omega_i,\sigma_i)$, choose a connected
component that contains a component of the vanishing locus and construct its 
moment map image as in the proof of Lemma~\ref{L:SimpleLoop}.
Recall that, by construction, the edge $e_2$ is vertical and contains
one fold point.
If the vertical line segment on the interior of the image with one
endpoint at the fold point of $e_2$ has its other vertex on an edge
of integral slope $m_i\in\bN$, then perform a symplectic cut, eliminating the 
fold point on edge $e_2$.
As in Example~\ref{E:ConnectedSum}, this is possible because the determinant
$\lvert \begin{smallmatrix} 0 & 1 \\ 1 & m_j \end{smallmatrix}\rvert$ has
norm $1$.

If not, proceed as follows, noting that on edges $e_p,\ldots e_n$ there are
no fold points.
\ben
\item Suppose there is an edge $e_j$, $j\ge p$, such that
\ben
\item $e_j$ has integral slope $m_j\ne 0$,
\item $e_j$ lies in a closed half-plane whose boundary contains the edge $e_2$,
and
\item the lower vertex of $e_j$ is higher than the fold point $p$.
\een
Then lengthen $e_j$ and $e_3$ (or $e_1$,
depending on which is on the opposite side of $l$ from $e_j$), 
scaling the edges
$e_4,\ldots, e_{j-1}$ (or $e_{j+1},\ldots e_n$) by a single constant so as to 
maintain an embedded polygonal path with folds.
\item 
\label{S:Deform}
Or, if there is an edge $e_j$, $j\ge p$, that has slope $0$ and lies above
$p$, then lengthen $e_{j+1}$ and $e_3$ (or $e_{j-1}$ and $e_1$,
depending on which side of $e_2$ the edge $e_j$ lies), scaling the edges
$e_4,\ldots, e_{j-1}$ (or $e_{j+1},\ldots e_n$) by a single constant so as to 
maintain an embedded polygonal path with folds.
\item Or, if there are no edges with integral slope, find the vertex with
maximal $y$ coordinate and name the edges on its left and right $e_j$ and 
$e_{j-1}$ with slopes $s_j=\frac{m_j}{n_j}$ and 
$s_{j-1}=\frac{m_{j-1}}{n_{j-1}}$. 
Remove a neighborhood of the corner, shortening $e_j$ and $e_{j-1}$ and
inserting a new edge $e'$ whose slope equals zero.  
If $m_j=1$ ($m_{j-1}=1$) then the boundary of the polygon is right polygonal
at the left (right) vertex of $e'$.  If not, the vertex defines an orbifold
singularity.  But any such singularity can be resolved equivariantly,
replacing the singular point with a union of spheres.
After resolving, the polygon will have
right polygonal boundary.  Next, go back to Step~\ref{S:Deform}.
\een
Now we may perform a symplectic cut to obtain $(X_{i+1},\omega_{i+1},
\sigma_{i+1}$). 
If we carry out this procedure $N=\abs{Z}$ times,
we obtain $(X_N,\omega_N,\sigma_N)$
for which $b_2^+(X_N)-\abs{Z_N}=c_N$.
Since this relation remains unchanged through all
of our constructions and $X$ is connected, $b_2^+(X_1)-\abs{Z_1}
= b_2^+(X)-\abs{Z_\omega}=1$.
\end{proof}

Proposition~\ref{P:SimplyConnectedT2} allows us to recognize
a simply connected $T^2$-manifold from its weighted orbit space, provided
the orbit space has at least five vertices.  
An ingredient in its proof is the calculation of the Euler characteristic
of a locally toric near-symplectic manifold from its weighted orbit space.

\begin{lemma}
\label{L:EulerChar}
Given a locally toric near-symplectic manfiold $(X,\omega,\sigma)$
with integral affine
base $(B,\cA,F)$, $\chi(X)=V$ where $V$ is the number of vertices
on the boundary of the integral affine base.
\end{lemma}
Recall that fold points do not count as vertices.  
\begin{proof}
  The total space of a locally toric manifold can be built up out of
  open sets, each of which is a neighborhood of a fiber.  This can be done 
  so that each neighborhood has Euler characteristic equal to $0$, 
  except for a 
  small neighborhood of the preimage of each vertex (which can be chosen to
  be a ball with Euler characteristic equal to $1$).  Furthermore, one can
  perform this operation so that as each neighborhood gets glued in, the
  gluing locus has Euler characteristic equal to $0$.
\end{proof}

\begin{proof}[Proof of Proposition~\ref{P:SimplyConnectedT2}]
If the orbit space has at least five vertices then it must be diffeomorphic
to a connected sum of copies of $\bC P^2$ and $\overline{\bC P}^2$
(cf.~\cite{OrlikRaymond}).
Then, invoking  Proposition~\ref{P:SimplestThmIII}
and Lemma~\ref{L:PosTurningMax}, we can calculate
\begin{align*}
m = b_2^+(X) & =\abs{Z_\omega}+1=\abs{F} +1 \\
             & = \frac{T-2\pi}{\pi}+1 = \frac{T}{\pi} -1.
\end{align*}
Then,
\begin{align*}
n = b_2^-(X) & =\chi(X) -b_2^+(X) -2 \\
             & = V-m-2.
\end{align*}
where the last equality follows from Lemma~\ref{L:EulerChar}.
\end{proof}

\subsection{$T^2$-manifolds with nontrivial fundamental groups}

The orbit space of any $T^2$ manifold whose fundamental group is nontrivial 
must also have nontrivial fundamental group.  In order to describe and
construct integral affine structures on orbit spaces of closed
manifolds we define
a few noncompact integral affine surfaces that serve as  building blocks.

\begin{defn}
Given an open interval $I$ in the positive reals, let $A_I$ be the annulus
$A_I = \{(x,y) | x^2 + y^2 \in I \}$. For any positive integer $q$,
define the {\em $q$-fold integral affine structure} $\cA_{I,q}$ on $A_I$ to 
be the pullback of the standard integral affine structure on $A_I \subset
\bR^2$ via the $q$-fold cover $(r,\theta) \mapsto (r,q\theta)$.
\end{defn}

\begin{defn}
  An integral affine surface with right polygonal boundary and edge 
  folds  $(B,\cA,F)$ is a {\em $q$-plug} if $B$ has exactly one end 
  modeled on the outer end of $A_{(a,b)}$ (for some interval 
  $(a,b)\subset \bR$), i.e. if there is an integral affine embedding 
  $\Phi:(A_{(a,b)},\cA_q)\ra (B\setminus F,\cA)$ such that 
  $B_0:=B\setminus \Phi(A_{(a,b)})$ is compact and any sequence of points
  in  $A_{(a,b)}$ converging to $r=a$ is sent via $\Phi$ to a 
  a sequence of points converging to $\partial B_0\subset B$.
\end{defn}

\begin{defn}
  An integral affine surface with right polygonal boundary and edge 
  folds  $(B,\cA,F)$ has a {\em $q$-hole} if
  $B$ has one end modeled on the inner end of $A_{(a,b)}$ (for some interval 
  $(a,b)\subset \bR$), i.e. if there is an integral affine embedding 
  $\Phi:(A_{(a,b)},\cA_q)\ra (B\setminus F,\cA)$ such that 
  $B_0:=B\setminus \Phi(A_{(a,b)})$ is connected and any sequence of points
  in  $A_{(a,b)}$ converging to $r=b$ is sent via $\Phi$ to 
  a sequence of points converging to $\partial B_0\subset B$.
\end{defn}

\begin{figure}
\begin{center}
\includegraphics[width=4in]{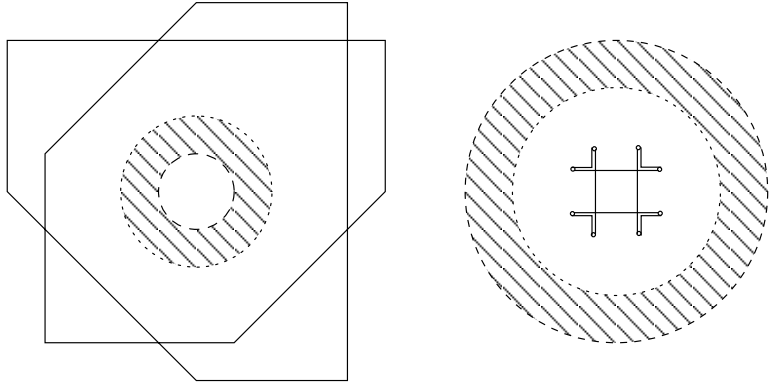}
\caption{On the left, an integral affine annulus with a $2$-hole, and
  on the right, a $3$-plug homeomorphic to an annulus. (In fact the
  figure on the right could also represent a $1$-plug homeomorphic to
  a twice-punctured torus.) }
\label{F:PlugAndHole}
\end{center}
\end{figure}

Figure~\ref{F:PlugAndHole} shows the immersed images of an integral
affine annulus with  a
$2$-hole (without any folds) and a $3$-plug (with four folds). The
shaded regions are the annular ends.

\begin{lemma} \label{L:Plugs}
Given any $q$ and any cyclic slope list $s_1,\ldots s_n$
that is right polygonal. there exists a $q$-plug
whose boundary realizes the given slope list.
\end{lemma}
\begin{proof}
  Without loss of generality, assume $s_1=0$ and add a slope, forming
  the slope list $s_1,s_2,\ldots s_n,s_{n+1}$ with $s_1=s_{n+1}=0$. Use
  Lemma~\ref{L:PolygonalPath} to construct a right polygonal path
  ending at $(x_0,y_0)$ and  
  representing this new slope list $s_1,s_2,\ldots s_n,s_{n+1}$.
  By construction, $(x_0,y_0)$ is in the interior of the first quadrant.
  Extend $e_2$ near the
  endpoint it shares with $e_3$ by the length $y_0$,
  thereby \lq\lq lowering\rq\rq\ the
  part of the path representing the slopes $s_3,\ldots s_n,s_{n+1}$.
  The result is a right polygonal path with endpoints at $(0,0)$ and 
  $(x_0,0)$.  Replace these two endpoints with fold points, concatenating
  the edges $e_{n+1}$ and $e_1$.    
  Although the resulting path is probably no longer embedded, it can be
  taken as the restriction to one boundary component of a continuous
  map of a closed annulus into $\bR^2$.  Furthermore, this continuous map can
  be assumed to be an immersion on the complement of a finite number of
  points  on that boundary component
  and to map the other boundary component to a large circle.  Then the
  interior of this annulus, with the integral affine structure induced
  by the immersion into $(\bR^2,\cA_0)$ is a $1$-plug.  To get a $q$-plug
  for $q > 1$, introduce $2(q-1)$ more folds in the edge $e_1$;
  see Figure~\ref{F:PlugConstruction}.
\begin{figure}
\begin{center}
\includegraphics[width=3in]{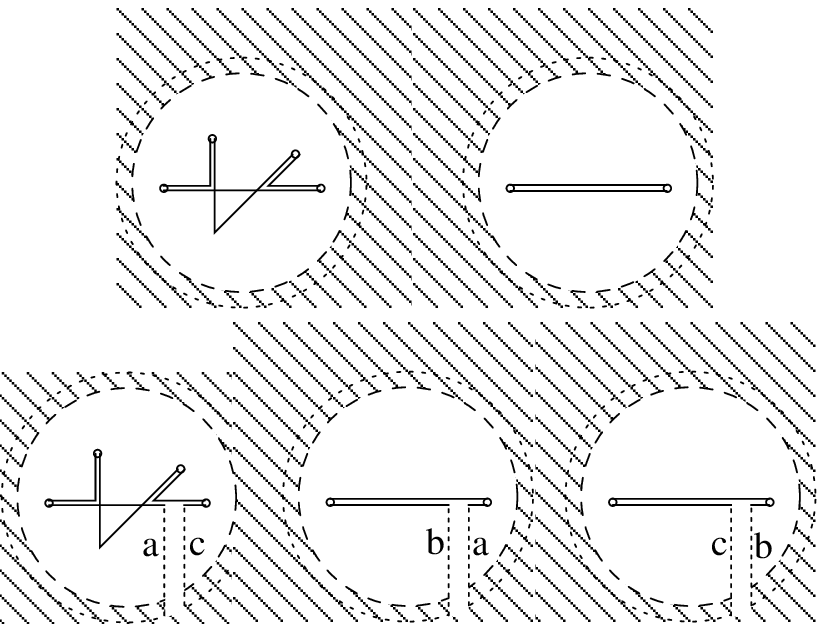}
\caption{Constructing a plug: Top left, a $1$-plug with three
  edges. Top right, a $1$-plug with one edge. On the bottom row, we
  indicate how to cut these open and glue them together (one of the
  first and two of the second) to make a $3$-plug with three edges
  (letters $a$, $b$, $c$ indicate gluing instructions).}
\label{F:PlugConstruction}
\end{center}
\end{figure}

\end{proof}

We now prove 
Theorem~\ref{T:HamiltActions} which asserts that there exists a 
near-symplectic structure making a torus action Hamiltonian only if
there is a component of the boundary of the weighted orbit space on
which the positive turning is at least $2\pi$; that this suffices if
the genus of the base is zero; and that a compact surface of any given 
genus and at least one boundary component 
can be the orbit space for a toric near-symplectic manifold. 

\begin{proof}[Proof of Theorem~\ref{T:HamiltActions}]

  Suppose $(X,\sigma)$ admits a near-symplectic structure with respect
  to which the action is Hamiltonian.  Let $(B,\cA,F)$ be the integral
  affine base and $\overline{\Phi}:B\ra\bR^2$ be the continuous
  extension of the integral affine immersion $\Phi:(B\setminus
  F,\cA)\ra(\bR^2,\cA_0)$ that defines the action on the symplectic
  locus.  The image under $\overline{\Phi}$ of each component of $\bdy
  B$ is the image of a polygonal loop with folds. Because $B$ is
  compact, its image $\overline{\Phi}(B)$ is a compact domain in
  $\bR^2$.  Consider an edge whose image contains points belonging to
  the boundary of $\overline{\Phi}(B)$.  Then that edge belongs to an
  immersed piecewise linear loop with folds $\gamma:[0,n]\ra\bR^2$
  that is the image of one component of $\bdy B$.
  
  For some indexing of the edges, there are 
  real numbers $a,b$ such that $[a,b]\subset[0,n]$ and $\gamma|_{[a,b]}$ 
  is an embedded loop.  
  Then $\gamma|_{[a,b]}$, with its orientation or the reverse, 
  bounds a disk inside the image of $\overline{\Phi}$.
  Because part of the image of $\gamma_{[a,b]}$ belongs to $\bdy\Phi(B)$,
  it must be  $\gamma|_{[a,b]}$ with its induced orientation that bounds
  this disk. 
  
  Let $\gint{a}$ 
  denote the greatest integer less than or equal to $a$ and $\lint{b}$
  denote the least integer greater than or equal to $b$.
  Then the positive turning of the path $\gamma|_{[\gint{a},\lint{b}]}$
  is greater than $\pi$.  The positive turning along the entire loop
  $\gamma$ must be at least as large, and an integer multiple of $\pi$.  
  Therefore,
  it is at least $2\pi$.  

  If $g=0$ then the base is a disk with $k$ holes for some $k\ge0$.  
  Choose a component of the boundary
  of the orbit space on which the total turning is at least $2\pi$.  Then
  use the proof of Lemma~\ref{L:SimpleLoop} to construct an embedded 
  right polygonal loop with folds representing the slope list for that
  boundary component.  This polygonal loop bounds a disk $D$ 
  in $\bR^2$.  Remove $k$ disjoint closed
  disks from the interior of $D$ so as to create $k$ $1$-holes.  
  Finally, glue in $k$ $1$-plugs, the boundaries of which
  realize the remaining cyclic slope lists encoded in the weighted orbit space.

  Meanwhile, the punctured torus example in Example~\ref{E:Assorted} 
  generalizes
  to give immersed examples for any $g$ and any $k \geq 1$, and hence
  to give examples of toric near-symplectic manifolds with these orbit spaces.
\end{proof}

A complete answer to the question of what $T^2$ manifolds admit
near-symplectic structures with respect to which the actions are 
Hamiltonians is not available, but is under investigation by the second
author.  However, in the next section we prove Theorem~\ref{T:locallyHamilt}
which provides a complete answer to an intermediate
question:  Given a $T^2$-action when does there exist a near-symplectic
form with respect to which the action is symplectic, and Hamiltonian
in a neighborhood of any orbit?

\section{Locally toric near-symplectic manifolds}
\label{S:LocallyToric}

As we have seen, an integral affine surface $(B,\cA,F)$
defines a toric near-symplectic manifold up to orbit preserving
symplectomorphism if and only if its boundary is right polygonal with folds
and there exists an integral affine immersion 
$\Phi:(B\setminus F,\cA)\ra (\bR^2,\cA_0)$.  

For the purposes of studying pseudoholomorphic curves in a symplectic 
$4$-manifold via $1$-complexes in a surface, the presence of a
global torus action is not necessary.  All one needs is 
 a {\em singular
Lagrangian fibration}\footnote{Loosely,
a singular Lagrangian fibration is a symplectic manifold $(X,\omega)$ 
together with a projection to a half-dimensional space $B$ such that 
over a dense open subset of $B$ the projection defines a locally trivial
fibration, each of whose fibers is Lagrangian.} in which the behavior of 
pseudoholomorphic curves
in the neighborhood of each singular fiber is understood.  
Locally toric fibrations comprise a convenient class of manifolds for this
purpose.  When the total space is symplectic, the list of manifolds that
admit locally toric fibrations is
short~\cite{LeungSymington.AlmostToric}. 
However, as we show in this section, there is a vast set of 
near-symplectic examples.

The following lemma implies that the base of a locally toric fibration 
is, like a toric fibration, equipped with a natural integral affine structure.
\begin{lemma}~\cite{Symington.4from2}
\label{L:FiberPreservingMaps}
Consider the Lagrangian fibration $\pi:(\bR^2\times T^2,dp\wedge dq)\ra\bR^2$
in which the map $\pi$ forgets the torus factor.
Suppose $U$
is a connected open subset of $\bR^2$.  Then an embedding $\Phi:
\mu^{-1}(U) \ra \bR^2\times T^2$ is a fiber-preserving symplectic
embedding if and only if it is of the form $(p,q) \mapsto (Ap + b,
A^{-T}q + f(p))$, where $A \in GL(2,\bZ)$, $b \in \bR^2$, $A^{-T}$ is
the inverse transpose of $A$, and $f: U \rightarrow T^2$ is a smooth
map such that $A^T \circ Df$ is symmetric.
\end{lemma}
Of course, the symplectic manifold $(\bR^2\times T^2,dp\wedge dq)$ also
supports the torus action $t\cdot(p,q)\ra (p,q+t)$ whose moment map is
$\pi$, so the equivariant symplectomorphic embeddings are
precisely those fiber-preserving embeddings with $A=Id$.

\begin{thm} \label{T:LocallyToricBases}
  Suppose $(B,\cA,F)$ is an
  integral affine surface with right polygonal boundary with edge folds.  
  Then $(B,\cA,F)$ is the integral affine base of a locally
  toric near-symplectic manifold.
\end{thm}

\begin{proof}
Cover $B$ by a union of open sets $\{U_\alpha\}$, 
each of which is contractible.
Then, for each $\alpha$ there exists an integral affine immersion
$\Phi_\alpha:(U_\alpha\setminus F,\cA)\ra(\bR^2,\cA_0)$ and hence 
a toric near-symplectic manifold $(X_\alpha,\omega_\alpha,\sigma_\alpha)$
with integral affine base $(U_\alpha,\cA,F)$,
constructed as in the proof of Theorem~\ref{T:ToricBases}.
By Lemma~\ref{L:AffineMoment}, on overlaps $U_\alpha\cap U_\beta$, the
two toric near-symplectic manifolds $(X_\alpha,\omega_\alpha,\sigma_\alpha)$
and $(X_\beta,\omega_\beta,\sigma_\beta)$ are orbit-preserving 
symplectomorphic (since $U_\alpha\cap U_\beta$ must be homeomorphic to
a union of contractible spaces).  
On triple intersections these gluing maps will be compatible, thereby
yielding a locally toric near-symplectic manifold.
\end{proof}

The only obstruction to the existence of a smooth global torus action
inducing a locally toric fibration is monodromy.
\begin{defn} \label{D:Monodromy} An integral affine structure $\cA$ on
  a surface $B$ determines a lattice $\Lambda(\cA)$ in $TB$ (coming
  via the defining atlas for $\cA$ from the standard integral lattice
  in $\bR^2$). The {\em monodromy} of $\cA$ is the monodromy 
  representation $\pi_1(B)  \ra GL(2,\bZ)$ of $\Lambda(\cA)$.
\end{defn}
Note that this is not the only obstruction if one requires the global
action to be Hamiltonian.  For example, consider
the square $\{(p_1,p_2) \,|\, |p_i|\le 1\}\subset (\bR^2,\cA_0)$
and identify the top and bottom edges to form an integral affine
cylinder.  The failure of this integral affine cylinder to isometrically
immerse in $(\bR^2,\cA_0)$ implies, by Theorem~\ref{T:ToricBases},
that it cannot be the integral affine base for a toric fibration.

In the locally toric setting, we do not always have uniqueness of the locally
toric manifolds with a given integral affine base. 
\begin{thm} \label{T:LocallyToricUniqueness}
An integral affine surface whose boundary is right polygonal with edge folds,
$(B,\cA,F)$, defines a unique locally toric near-symplectic manifold
if and only if $B$ is either noncompact or has nonempty boundary (i.e.
has the homotopy type of a $1$-complex).  
The uniqueness is up to fiber-preserving
homeomorphism that is a  symplectomorphism on the symplectic locus.
\end{thm}
Note that the base $B$ need not be orientable.

\begin{proof}
If $B$ has the homotopy type of a $1$-complex then we can find in $B$
a collection of disjoint properly embedded arcs $\{\gamma_\alpha\}$ in
$B\setminus F$ such that $B\setminus \cup_\alpha \gamma_\alpha$ is a
disjoint union of simply connected surfaces.  Choose disjoint open
collar neighborhoods $\{V_\alpha\}$ of the $\{\gamma_\alpha\}$ such
that for each $\alpha$, $V_\alpha\subset B\setminus F$.  Define open
sets $\{U_\beta\}$ such that each $U_\beta$ is the union of one
component of $B\setminus \cup_\alpha \gamma_\alpha$ and all of the
$V_\alpha$ that have nonempty intersection with that
component. Arrange that each $U_\beta$ is simply connected by going
back and including more arcs in the set $\{\gamma_\alpha\}$ if
necessary.

By Lemma~\ref{L:SimplyConnected}, each of these integral affine
surfaces $(U_\beta\setminus F,\cA)$ immerses isometrically in
$(\bR^2,\cA_0)$ and hence defines a unique toric near-symplectic
manifold, say $(X_\beta,\omega_\beta,\sigma_\beta)$.  Then all of the
locally toric near-symplectic manifolds defined by $(B,\cA,F)$ can be
built out of the $(X_\beta,\omega_\beta,\sigma_\beta)$ by gluing maps
between neighborhoods that project to the $V_\alpha$.  By
Lemma~\ref{L:FiberPreservingMaps}, each of these maps can be
expressed, in local coordinates on the top dimensional fibers, as
$(p,q)\mapsto (Ap+b,A^{-T}q+\phi(p))$ for some $(A,b)\in{\rm
Aff}(2,\bZ)$ and some $\phi(p)$ that is the time-one flow of a
Hamiltonian vector field.  The element $(A,b)\in {\rm Aff}(2,\bZ)$ is
uniquely determined by the integral affine structure on $B$, while the
arguments in the proof of Theorem~\ref{T:ToricUniqueness} show that
particular choices of $\phi$ have no effect on the global structure.
Thus the global structure defined by $(B,\cA,F)$ is unique.

The hypothesis that $B$ have the homotopy type of a $1$-complex is
necessary because there exist closed integral affine surfaces that
each are the integral affine base of more than one locally toric
manifold.  This is true even if $H^2(B,\bZ)=0$ as evidenced by the
existence of multiple locally toric manifolds whose integral affine
base is a single integral affine Klein bottle, as shown
in~\cite{LeungSymington.AlmostToric}.
\end{proof}

We now take up the problem of realizing 
$T^2$-manifolds as locally toric with respect to some near symplectic
structure, i.e. constructing integral affine surfaces with right
polygonal boundary that define the underlying fibration.

\begin{lemma} \label{L:GenusAndHoles} For any $g \geq 1$ and $k \geq
  1$, there exists a noncompact integral affine surface $B$ with
  empty boundary, genus $g$ and $k$ ends, each of which 
  is a $q$-hole for some $q$ (not necessarily the same
  $q$ for each end).
\end{lemma}

\begin{proof}
  For an appropriately chosen positive integer $p$, let $R$ be a
  $p$-by-$1$ rectangle in $\bR^2$ with quarter- and half-circles, 
  all of radius $r<\frac{1}{2}$,
  removed at the corners of the $p$ $1$-by-$1$ rectangles making up
  $R$. Use $R$ as a fundamental domain to build $B$, gluing
  appropriate edges to each other via translations. (Without removing
  the quarter- and half-circles  this is a method to produce a flat
  metric on a closed genus $g$ surface with $k$ singular points; the
  first author learned this trick from A. Abrams, who claims to have
  learned it from~\cite{Thurston.Surface}.)
Note that, to avoid monodromy, the left edge must be glued to the right, 
and edges on the top must be glued to edges on the bottom.  
Figure~\ref{F:GenusAndHoles} shows an example for $g=2$,
$k=2$, with one end having $q=1$ and the other having $q=3$.
\begin{figure}
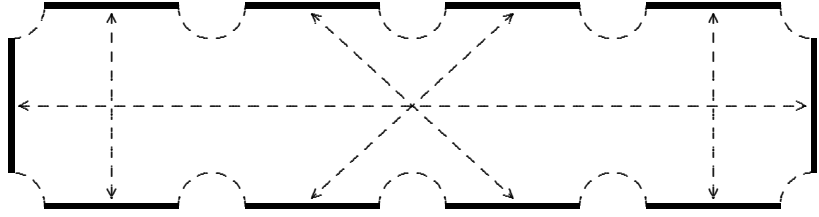

\begin{center}
\include{GenusAndHoles}
\caption{A fundamental domain for a genus $2$ surface with holes; arrows
  indicate gluing rules.}
\label{F:GenusAndHoles}
\end{center}
\end{figure}
\end{proof}

\begin{lemma} \label{L:ParallelAnnulus} Given any slope $s\in\bQ\cup\infty$, 
there exists an integral affine annulus with trivial monodromy, with each
  boundary consisting of a single edge with slope $s$.
\end{lemma}

\begin{proof}
Glue two sides of a parallelogram via a translation.
\end{proof}

\begin{lemma} \label{L:NonParallelAnnulus} Given two right polygonal cyclic
  slope lists, one of which has at least two slopes, 
  there exists an integral affine annulus
  with trivial monodromy whose boundary realizes the two slope lists.
\end{lemma}

\begin{proof}
Let the slope lists be $s_1,\ldots, s_k$ and $t_1,\ldots, t_l$.
If one list has only one slope, then let $\tau$ be a translation in
the direction of this slope such that $(x_0,y_0):=\tau((0,0))$ is in the closed
right half plane.
 Otherwise choose such a translation $\tau$ arbitrarily.
If $x_0,y_0>0$ then construct, 
as in the proof of Lemma~\ref{L:PolygonalPath}, two right polygonal
paths with folds representing $s_1,\ldots, s_k,s_1$ and $t_1,\ldots, t_l, t_1$
starting at $(0,0)$ and ending at $(x_0,y_0)=\tau((0,0))$.
If $x_0>0$ and $y_0\le 0$, then construct the path with endpont $(x_0,y_1)$
for some $y_1>0$ and lengthen the edge $e_2$ so as 
to lower the endpoint from $y_1$ to $y_0$.  If $x_0=0$, then interchange the
roles of $x$ and $y$ (thereby affecting the slopes also), 
make the construction, and switch back.

Construct two $\tau$-invariant periodic paths by concatenating tranlated
copies of these paths.
Translate one of the periodic paths so they become disjoint, and rotate one by
$180^\circ$ so that they are the oriented boundary of the 
strip $S$ in between. Then $S/\langle\tau\rangle$ is a base
with folds homeomorphic to an annulus with boundary realizing the two
slope lists.
\end{proof}

\begin{lemma} \label{L:PuncturedAnnulus} Given three slopes $s_1, s_2,
  s_3 \in \bQ \cup \{\infty\}$, there exists an integral affine
  twice-punctured disk with trivial monodromy with each boundary
  component consisting of a single edge of slope $s_i$.
\end{lemma}

\begin{proof}
If two of the slopes are equal, use Lemma~\ref{L:ParallelAnnulus} with
one $1$-hole, and fill with a plug from Lemma~\ref{L:Plugs}.
Otherwise, without loss of generality we can assume that $s_1 < 0 =
s_2 < s_3 < \infty$  (recall that the order is
unimportant). Figure~\ref{F:3Slopes} then illustrates the construction
by means of a fundamental domain embedded in $\bR^2$.

\begin{figure}
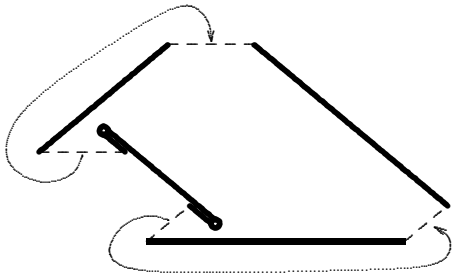

\begin{center}
\include{3SlopesV2}
\caption{Building a twice-punctured disk with one edge per boundary.}
\label{F:3Slopes}
\end{center}
\end{figure}

\end{proof}

\begin{proof}[Proof of Theorem~\ref{T:locallyHamilt}]
Theorems~\ref{T:LocallyToricBases} 
and~\ref{T:LocallyToricUniqueness} reduce our proof to proving that,
given a genus $g$ and $k$ slope lists, there exists an 
immersion of a genus $g$ surface with $k$ boundary components -- or a
fundamental domain of such a surface -- 
into $\bR^2$ so that its boundary is a union of right polygonal 
loops with folds realizing the prescribed slope lists if and only if the
hypotheses of the theorem are satisfied.

Notice that since the union of principal orbits is an oriented $T^2$-bundle, 
the integral affine structure must have trivial monodromy.
Accordingly, the gluing maps used to construct the surface from an
immersion of a fundamental domain must be translations.

We now go through the various cases of specific genera and numbers of 
boundary components.

If the surface has no boundary, then the presence of an integral affine 
structure implies that the tangent bundle admits a flat connection, and 
therefore has zero Euler characteristic.  Since the monodromy is trivial,
the surface must be orientable, and therefore a torus.
Any torus bundle over a torus with trivial monodromy (but not necessarily 
trivial Chern class) supports a near-symplectic structure with respect to 
which the action is locally toric.  Since the boundary is empty, the 
near-symplectic structure is in fact symplectic.  Such manifolds were 
classified by Mishachev~\cite{Mishachev}. The bases of such locally toric 
symplectic manifolds can be constructed by choosing a parallelogram in the 
plane that is integral affine equivalent to a recangle, and pairwise 
identifing the opposite edges via translations.

Now suppose the surface has nonempty boundary.  We first suppose that $g\ge 1$.
Then Lemma~\ref{L:GenusAndHoles} tells us how to construct a noncompact
integral affine surface (without boundary) with $k$ ends, and
Lemma~\ref{L:Plugs} asserts that we can glue into each of these ends the
collar neighborhood of a boundary component that realizes any given slope
list.

If the surface has genus zero we have three cases to consider, $k\ge 3$,
$k=2$, and $k=1$.  

Suppose $k\ge 3$.  If all of the slope lists consist of exactly one
slope, we can use Lemma~\ref{L:PuncturedAnnulus} and any three of the
slopes to construct a punctured annulus realizing these three slopes,
and then remove neighborhoods of slits and use Lemma~\ref{L:Plugs} to
fill in the remaining boundary components.  Otherwise,
Lemma~\ref{L:NonParallelAnnulus} tells us how to construct an annulus
realizing two of the slope lists, one of which has at least two
slopes; then again, we can remove disks to create $1$-holes and then
glue in $1$-plugs realizing the remaining slope lists.

If $k=2$ and one slope list has at at least two slopes (which happens if and
only if $T_0 \geq \pi$) then we can use  Lemma~\ref{L:NonParallelAnnulus} 
together with $1$-plugs from Lemma~\ref{L:Plugs} to construct the required
immersed surface.  Otherwise, if both slope lists consist of just one slope, 
then we merely need to invoke Lemma~\ref{L:ParallelAnnulus}.  In this
latter case $T_0=0$.  Notice that the two slopes must be equal for otherwise
we would be trying to glue the two parallel sides of a trapezoid via a
translation -- which we can do only if they have the same length, i.e. if
the trapezoid is actually a rectangle.

For the last case, suppose $k=1$.  In this case the weighted orbit space
is simply connected.  Therefore given any integral affine structure 
$\cA$ on the orbit space, there is an integral affine immersion 
$(B,\cA)\ra(\bR^2,\cA_0)$.  Consequently, if the action is locally Hamiltonian
it must be Hamiltonian.  Therefore, the necessary and sufficient
conditions are contained in Theorem\ref{T:HamiltActions}.

Lastly, to compute the number of components of the vanishing locus,
note that $|Z_\omega| = \abs{F}$, the number of folds. The total turning on
all boundaries of the base is $T-\abs{F}\pi$, where $T$ is the sum of the
positive turnings over all boundary components.  The Gauss-Bonnet Theorem
implies that $2\pi\chi=T-\pi \abs{F}$, and hence 
$|Z_\omega| = \frac{1}{\pi}T-2\chi)$.
\end{proof}

We now forget about torus actions. 
In order to highlight the plenitude of closed locally toric near-symplectic
manifolds we prove Propositions~\ref{P:AnyMonodromy} and 
\ref{P:InfiniteFamily}.  Recall that the first proposition asserts
that any monodromy representation of
the fundamental group of any surface can arise as the monodromy of a
locally toric fibration of a closed near-symplectic manifold, while the
second asserts that there is an infinite family of mutually non-diffeomorphic
closed near-symplectic manifolds that support locally toric fibrations having
nontrivial monodromy, none of
which could support a locally toric fibration with trivial monodromy.

\begin{proof}[Proof of Proposition~\ref{P:AnyMonodromy}]
  Suppose the free group has $n$ generators $x_1, \ldots, x_n$, and
  the homomorphism maps $x_i$ to the matrix $(\begin{smallmatrix} a_i
    & c_i \\ b_i & d_i
\end{smallmatrix} ) \in GL(2,\bZ)$. 
Draw a (probably nonconvex) polygon with at least $3n$ edges in
$\bR^2$ having rational slope, such that at every vertex the two
incident edges have primitive integral tangent vectors $v,w$
satisfying $\det(v,w) = \pm 1$. Draw the polygon so that it has $2n$
distinguished edges $e_1, \ldots, e_n$ and $f_1, \ldots, f_n$,
satisfying the following properties:
\begin{enumerate}
\item The primitive integral tangent vector to each $e_i$ is
$(\begin{smallmatrix} 1 \\ 0 \end{smallmatrix} )$, while the
primitive integral tangent vector to the edge immediately preceding
$e_i$ (with the boundary orientation) is $(\begin{smallmatrix} 0 \\
-1 \end{smallmatrix} )$ and the primitive integral tangent vector to
the edge immediately following $e_i$ (with the boundary orientation)
is $(\begin{smallmatrix} 0 \\ 1 \end{smallmatrix} )$.
\item The primitive integral tangent vector to each $f_i$ is
$(\begin{smallmatrix} a_i \\ b_i \end{smallmatrix} )$, while the
primitive integral tangent vector to the edge immediately preceding
$f_i$ (with the boundary orientation) is $(\begin{smallmatrix} c_i \\
d_i \end{smallmatrix} )$ and the
primitive integral tangent vector to the edge immediately following
$f_i$ (with the boundary orientation) is $(\begin{smallmatrix} -c_i \\
-d_i \end{smallmatrix} )$.
\item The lengths of the $e_i$ and $f_i$ are chosen so that, for each $i$, 
the Euclidean lengths of $(\begin{smallmatrix} a_i
    & c_i \\ b_i & d_i
\end{smallmatrix} ) e_i$ and $f_i$ are equal.
\item As one traverses the boundary of the polygon counterclockwise, the
edges $e_i$ and $f_i$ are ordered so that identification of each $e_i$
with each $f_i$ produces the surface $B$.
\end{enumerate}

Now construct a base with folds from this polygon by replacing each
concave corner with a convex corner and a nearby fold (extend one
incident edge a little past the corner and then immediately double
back, creating a fold point, then turn onto the other incident edge)
and by gluing a rectangular neighborhood of $e_i$ to a rectangular
neighborhood of $f_i$ via $(\begin{smallmatrix}a_i & c_i \\ b_i & d_i
\end{smallmatrix} )$ followed by a suitable translation.
\end{proof}

\begin{proof}[Proof of Proposition~\ref{P:InfiniteFamily}]
  Figure~\ref{F:NonTrivialFamily} gives the construction of the
  infinite family, where the ellipses in the middle are to be
  interpreted as representing $n$ \lq\lq slits.\rq\rq\ 
  Call these manifolds $X_n$, $n\ge 0$.  

  These manifolds are mutually non-diffeomorphic because adding 
  a slit has the effect of adding a generator and no relations to the
  fundamental group of the $4$-manifold.
  
  To see why no $X_n$ can be diffeomorphic to a locally toric
near-symplectic manifold whose fibration has trivial monodromy we note
that, by Lemma~\ref{L:EulerChar}, $\chi(X_n)=1$ for all $n$.  However,
a near-symplectic manifold equipped with a locally toric fibration
having trivial monodromy can never have Euler characteristic equal to
$1$.  Indeed, this would imply that the integral affine base would
have a boundary component with just one vertex.  With trivial
monodromy, this means that the isotropy subgroups for orbits whose
images are on either side of the vertex are the same, violating the
requirement that the boundary of the integral affine base have right
polygonal boundary.
\begin{figure}
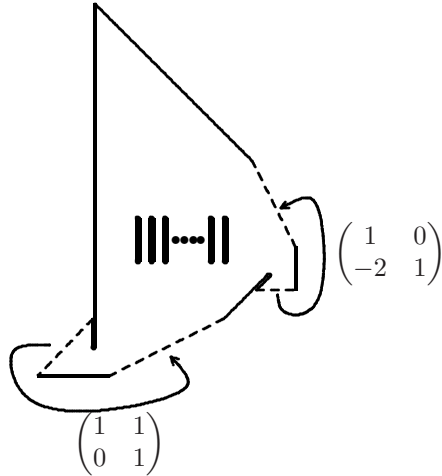

\begin{center}
\include{NonTrivialFamily}
\caption{An infinite family of examples with nontrivial
  monodromy. (All lines in this figure have slope $0$, $\pm 1/2$, $\pm
  1$, $\pm 2$ or $\infty$.)}
\label{F:NonTrivialFamily}
\end{center}
\end{figure}

\end{proof}

\appendix

\section{Hamiltonian actions}
\label{A:Hamiltonian}

In this section we prove Proposition~\ref{P:Hamiltonian} which gives a 
simplified definition of a Hamiltonian action when the group is a torus
(what we called a \lq\lq topologist's definition\rq\rq).

In general, given a group action $\sigma:G\times X\ra X$ on a symplectic
manifold $(X,\omega)$, there are two types of natural vector fields.  For
each $\xi$ in the Lie algebra $\cG$ there is the infinitesimal action
$V_\xi$, while for each smooth function $f:X\ra\bR$ there is the
Hamiltonian vector field $V_f$ defined by $\omega(V_f,\cdot)=-df$.  

\begin{defn}\label{D:Hamiltonian}
A group action $\sigma:G\times X\ra X$ on a symplectic manifold
$(X,\omega)$ is {\em Hamiltonian} if there is a Lie algebra
homomorphism from $\cG$ to $C^\infty(X)$ (equipped with the Poisson
bracket) that sends $\xi$ to the function $ f_\xi$ so that
$V_{f_\xi}=V_\xi$.
\end{defn}

\begin{proof}[Proof of Proposition~\ref{P:Hamiltonian}]
When the group is a torus $T^n$, $\cG=\cG^*=\bR^n$.  In that case,
Definition~\ref{D:Hamiltonian} implies the existence of a
moment map $\mu:M\ra\cG^*$ defined implicitly by $\langle \mu(x),\xi
\rangle =f_\xi(x)$. Since $df_\xi(W)=\xi\cdot D\mu(W)$ for any vector
$W$, one direction of the proposition is clear.

For the converse, the functions $f_\xi$ are just $\xi\cdot \mu$.  We
only need to check that $\xi\mapsto f_\xi$ is a Lie algebra
homomorphism.  Since the action is abelian, for any
$\xi,\eta\in\bR^2=\cT$ we have $[\xi,\eta]=0$.  Therefore, $\xi\mapsto
f_\xi$ will be a Lie algebra homomorphism provided
$\{f_\xi,f_\eta\}=\omega(V_\xi,V_\eta)=0$.

First of all, for any $\xi,\eta$, $\omega(V_\xi,V_\eta)$ is constant
on the orbit of $V_\eta$.  This follows from two facts: the time $t$
flow $\phi_t$ of $V_\eta$ is a symplectomorphism for each $t\in\bR$
and the action is abelian so $(\phi_t)_*V_\nu=V_\nu\circ\phi_t$ for
each $\nu\in\cT$.  Calculating $\omega(V_\xi,V_\eta) = -df_\xi(V_\eta)
= -V_\eta f_\xi$ we see that $V_\eta f_\xi$ is constant on the orbit
of $V_\xi$.  For almost all $\eta$, the orbit is a circle, and hence
the derivative $V_\eta f_\xi$, being equal to a constant, must be equal
to zero, implying that $f_\xi$ is constant on the orbit of $V_\eta$.
Continuity of the moment map then implies that this is true for all
$\eta$, thereby establishing that $\omega(V_\xi,V_\eta)=0$ for all
$\xi,\eta$.
\end{proof}

\begin{rmk}
The compactness of the torus is essential here.  Note that there is a
symplectic action of $\bR^2$ on $(\bR^2,\omega_0)$, by translation,
that satisfies the \lq\lq topologist's definition\rq\rq\
(Proposition~\ref{P:Hamiltonian}) but is certainly not Hamiltonian.
\end{rmk}

\bibliographystyle{gtart}
\bibliography{bibliography}

\end{document}